\documentclass{amsart} 
\usepackage{amsmath,amsfonts,amssymb,amscd,mathrsfs}
\usepackage{graphicx,comment,alltt}
\usepackage{mathptmx,courier}
\usepackage[scaled=0.95]{helvet}

\specialcomment{Notes}{\begin{quote}\footnotesize\begin{notes}\end{notes}}{\end{quote}}
\usepackage{ifpdf}
\ifpdf
\usepackage[%pdftex,
            pdftitle={Knot colouring polynomials},
            pdfauthor={Michael Eisermann},
            pdfstartview=FitV,
            bookmarksopen=true,
            colorlinks=true, 
            linkcolor=blue, 
            citecolor=blue, 
            urlcolor=blue,
            ]{hyperref}
\else
\usepackage[dvips]{hyperref}
\fi

\title{Knot colouring polynomials}

\author{Michael Eisermann}
\address{Institut Fourier, Universit\'e Grenoble I, France}
\email{Michael.Eisermann@ujf-grenoble.fr}
\urladdr{www-fourier.ujf-grenoble.fr/{\textasciitilde}eiserm}

\date{\today}

\theoremstyle{plain}
\newtheorem{theorem}{Theorem}[section]
\newtheorem{lemma}[theorem]{Lemma}
\newtheorem{proposition}[theorem]{Proposition}
\newtheorem{corollary}[theorem]{Corollary}

\theoremstyle{definition}
\newtheorem{definition}[theorem]{Definition}
\newtheorem{remark}[theorem]{Remark}
\newtheorem{example}[theorem]{Example}
\newtheorem{question}[theorem]{Question}

\newcommand{\N}{\mathbb{N}}                       % natural numbers
\newcommand{\Z}{\mathbb{Z}}                       % integers
\newcommand{\R}{\mathbb{R}}                       % real numbers
\newcommand{\K}{\mathbb{K}}                       % ground ring
\newcommand{\A}{\mathbb{A}}                       % ring extension 

\newcommand{\Hom}{\operatorname{Hom}}             % set of homomorphisms
\newcommand{\compose}{\circ}                      % composition of maps
\newcommand{\tensor}[1][]{\otimes_{#1}}           % tensor product [over a ring]
\newcommand{\id}{\operatorname{id}}               % identity map
\newcommand{\tr}{\operatorname{tr}}               % trace of an endomorphism
\newcommand{\End}{\operatorname{End}}             % set of endomorphisms

\newcommand{\lto}[1][]{\mathrel{\smash{\overset{\smash{#1}}{\longrightarrow}}}}
\newcommand{\conj}{\operatorname{conj}}           % conjugation

\renewcommand{\S}{\mathbb{S}}                     % unit sphere
\newcommand{\D}{\mathbb{D}}                       % closed unit ball
\newcommand{\minus}{\smallsetminus}               % set complement

\newcommand{\knots}{\mathscr{K}}                  % set of knots
\newcommand{\links}{\mathscr{L}}                  % set of knots
\newcommand{\lk}{\operatorname{lk}}               % linking number
\newcommand{\uk}{\bigcirc}                        % the unknot
\newcommand{\cs}{\mathbin{\sharp}}                % connected sum
\newcommand{\rev}{{!}}                            % reverse knot
\newcommand{\obv}{{\times}}                       % obverse knot
\newcommand{\inv}{{\ast}}                         % inverse knot

\newcommand{\TF}[1]{\smash{F_{\!#1}}}             % total colouring number
\newcommand{\F}[2]{\smash{F^{#1}_{\!#2}}}         % colouring number
\renewcommand{\P}[2]{\smash{P^{#1}_{\!#2}}}       % colouring polynomial
\newcommand{\Col}{\operatorname{Col}}             % set of quandle colourings
\newcommand{\tsa}{\mathbin{\overline{\ast}}}      % inverse quandle operation
\newcommand{\weight}[2]{\langle#1|#2\rangle}      % weight of a colouring
\newcommand{\CSSI}[2]{\smash{\mathit{S_{#1}^{#2}}}}% cocycle state-sum invariant
\newcommand{\Ext}{\mathscr{E}}                    % equivalence classes of extensions

\newcommand{\GF}[1][]{\mathbb{F}_{#1}}            % finite field
\newcommand{\Gal}{\operatorname{Gal}}             % galois group
\newcommand{\Br}[1]{\operatorname{B}_{#1}}        % braid group
\newcommand{\Sym}[1]{\operatorname{S}_{#1}}       % symmetric group
\newcommand{\Alt}[1]{\operatorname{A}_{#1}}       % alternating group
\newcommand{\PSL}[2]{\operatorname{PSL}_{#1}{#2}} % projective special linear group
\newcommand{\Aut}{\operatorname{Aut}}             % set of automorphisms
\newcommand{\Inn}{\operatorname{Inn}}             % set of inner automorphisms
\newcommand{\gen}[1]{\langle{#1}\rangle}          % subgroup given by generators

\begin{document} %%%%%%%%%%%%%%%%%%%%%%%%%%%%%%%%%%%%%%%%%%%%%%%%%%%%%%%%%%%%
%%%%%%%%%%%%%%%%%%%%%%%%%%%%%%%%%%%%%%%%%%%%%%%%%%%%%%%%%%%%%%%%%%%%%%%%%%%%%

\begin{abstract}
  This article introduces a natural extension of colouring numbers of knots,
  called colouring polynomials, and studies their relationship to Yang-Baxter
  invariants and quandle $2$-cocycle invariants.

  For a knot $K$ in the $3$-sphere let $\pi_{K}$ be 
  the fundamental group of the knot complement $\S^3\minus K$, 
  and let $m_K,l_K \in \pi_{K}$ be a meridian-longitude pair.
  Given a finite group $G$ and an element $x\in G$ we consider
  the set of representations $\rho\colon \pi_{K} \to G$ with $\rho(m_K) = x$
  and define the colouring polynomial $\P{x}{G}(K) := \sum_\rho \rho(l_K)$.
  The resulting invariant maps knots to the group ring $\Z G$.
  It is multiplicative with respect to connected sum and 
  equivariant with respect to symmetry operations of knots.
  Examples are given to show that colouring polynomials 
  distinguish knots for which other invariants fail, in particular 
  they can distinguish knots from their mutants, obverses, inverses, or reverses.
  
  We prove that every quandle $2$-cocycle state-sum invariant of knots 
  is a specialization of some knot colouring polynomial.  This provides 
  a complete topological interpretation of these invariants 
  in terms of the knot group and its peripheral system.
  Furthermore, we show that $\P{x}{G}$ can be presented as a Yang-Baxter 
  invariant, i.e.\ as the trace of some linear braid group representation.
  This entails in particular that Yang-Baxter invariants 
  \emph{can} detect non-inversible and non-reversible knots.
\end{abstract}

 % Invariants of knots and 3-manifolds


%%%%%%%%%%%%%%%%%%%%%%%%%%%%%%%%%%%%%%%%%%%%%%%%%%%%%%%%%%%%%%%%%%%%%%%%%%%%%

% \headline{20mm}{\eprintinfo}
\vspace*{-10mm}

\maketitle

\vspace*{-2mm}

%%%%%%%%%%%%%%%%%%%%%%%%%%%%%%%%%%%%%%%%%%%%%%%%%%%%%%%%%%%%%%%%%%%%%%%%%%%%%

\section{Introduction and statement of results} \label{sec:Introduction}

To each knot $K$ in the $3$-sphere $\S^3$ we can associate its knot group, 
that is, the fundamental group of the knot complement, 
denoted by $\pi_{K} := \pi_1(\S^3\minus K)$.
This group is already a very strong invariant: it classifies 
unoriented prime knots \cite{Whitten:1987,GordonLuecke:1989}.
In order to capture the complete information, 
we consider a meridian-longitude pair $m_K,l_K \in \pi_{K}$:
the group system $(\pi_{K},m_K,l_K)$ classifies 
oriented knots in the $3$-sphere \cite{Waldhausen:1968}.
In particular, the group system allows us to tackle
the problem of detecting asymmetries of a given knot
(see \textsection\ref{sub:KnotGroupSymmetries}).
Using this ansatz, M.\,Dehn \cite{Dehn:1914} proved in 1914 that 
the two trefoil knots are chiral, and, half a century later, 
H.F.\,Trotter \cite{Trotter:1963} proved % in 1963 
that bretzel knots are non-reversible. 
We will recover these results using knot colouring 
polynomials (see \textsection\ref{sub:Examples}).

Given a knot $K$, say represented by some planar diagram $D$, 
we can easily read off the Wirtinger presentation of $\pi_{K}$ 
in terms of generators and relations (see \textsection\ref{sub:WirtingerQuandles}). 
In general, however, such presentations are very difficult to analyze.
As R.H.\,Crowell and R.H.\,Fox \cite[\textsection VI.5]{CrowellFox:1963} put it:
\begin{quote}
  ``What is needed are some standard procedures 
  for deriving from a group presentation some 
  easily calculable algebraic quantities 
  which are the same for isomorphic groups 
  and hence are so-called group invariants.''
  % (R.H.\,Crowell and R.H.\,Fox \cite[\textsection VI.5]{CrowellFox:1963})
\end{quote}
The classical approach is, of course, to consider 
abelian invariants, most notably the Alexander polynomial.
In order to effectively extract non-abelian information,
we consider the set of knot group homomorphisms 
$\Hom(\pi_{K};G)$ to some finite group $G$. 
The aim of this article is to organize this information
and to generalize colouring numbers to colouring polynomials.
In doing so, we will highlight the close relationship
to Yang-Baxter invariants and their deformations on the one hand,
and to quandle cohomology and associated state-sum invariants on the other hand.

\subsection{From colouring numbers to colouring polynomials}

A first and rather crude invariant is given by 
the total number of $G$-representations, denoted by
\[
\TF{G}(K) := |\Hom(\pi_{K};G)|.
\]
This defines a map $\TF{G} \colon \knots \to \Z$
on the set $\knots$ of isotopy classes of knots in $\S^3$.
This invariant can be refined by further specifying 
the image of the meridian $m_K$, that is,
we choose an element $x\in G$ and consider only
those homomorphisms $\rho\colon \pi_{K} \to G$ satisfying $\rho(m_K) = x$.
Their total number defines the knot invariant 
\[
\F{x}{G}(K) := |\Hom(\pi_{K},m_K;G,x)|.
\]

\begin{example}
  Let $G$ be the dihedral group of order $2p$, where $p\ge3$ is odd,
  and let $x\in G$ be a reflection. Then $\F{x}{G}$ is the number of $p$-colourings
  as introduced by R.H.\,Fox \cite{Fox:1962}, here divided by $p$ 
  for normalization such that $\F{x}{G}(\uk) = 1$.
\end{example}

We will call $\F{x}{G}$ the \emph{colouring number} associated with $(G,x)$,
in the dihedral case just as well as in the general case of an arbitrary group.
Obviously $\TF{G}$ can be recovered from $\F{x}{G}$ by summation over all $x\in G$.
In order to exploit the information of meridian \emph{and} 
longitude, we introduce knot colouring polynomials as follows:

\begin{definition} \label{def:ColPoly}
  Suppose that $G$ is a finite group and $x$ is one of its elements.
  The \emph{colouring polynomial} $\P{x}{G}\colon \knots \to \Z{G}$ 
  is defined as 
  \[
  \P{x}{G}(K) := \sum_\rho \; \rho(l_K) ,
  \]
  where the sum is taken over all homomorphisms 
  $\rho\colon \pi_{K} \to G$ with $\rho(m_K) = x$. 
\end{definition}

By definition $\P{x}{G}$ takes its values in the semiring $\N{G}$,
but we prefer the more familiar group ring $\Z{G} \supset \N{G}$.
Note that we recover the colouring number $\F{x}{G} = \varepsilon \P{x}{G}$
by composing with the augmentation map $\varepsilon\colon \Z{G} \to \Z$.
As it turns out, colouring polynomials allow us in a simple 
and direct manner to distinguish knots from their mirror images, 
as well as from their reverse or inverse knots.
We will highlight some examples below.

\subsection{Elementary properties}

The invariant $\P{x}{G}$ behaves very much like classical knot polynomials.
Most notably, it nicely reflects the natural operations on knots:
$\P{x}{G}$ is multiplicative under connected sum and equivariant 
under symmetry operations (\textsection\ref{sub:KnotGroupSymmetries}). 

Strictly speaking, $\P{x}{G}(K)$ is, of course, 
not a polynomial but an element in the group ring $\Z{G}$.
Since $l_K$ lies in the commutator subgroup $\pi_{K}'$ and commutes with $m_K$,
possible longitude images lie in the subgroup $\Lambda = C(x) \cap G'$.
Very often this subgroup will be cyclic, $\Lambda = \gen{t}$ say, in which case $\P{x}{G}$
takes values in the truncated polynomial ring $\Z\Lambda = \Z[t]/(t^n)$.
Here is a first and very simple example:

\begin{example} \label{exm:TrefoilAlt5}
  We choose the alternating group $G = \Alt{5}$ with basepoint $x = (12345)$.
  Here the longitude subgroup $\Lambda = \gen{x}$ is cyclic of order $5$.
  The colouring polynomials of the left- and right-handed trefoil knots 
  are $1+5x$ and $1+5x^{-1}$ respectively, hence the trefoil knots are chiral.
  (A typical colouring is shown in \textsection\ref{sec:Quandles}, 
  Figure \ref{fig:TrefoilAlt5}.)

  Starting from scratch, i.e.\ knot diagrams and Reidemeister moves,
  one usually appreciates Fox' notion of $3$-colourability \cite{Fox:1962}
  as the simplest proof of knottedness.  In this vein, the preceding example 
  is arguably one of the most elementary proofs of chirality,
  only rivalled by Kauffman's bracket leading to 
  the Jones polynomial \cite{Kauffman:1987}.
  % and perhaps the calculation of the signature from 
  % the Seifert matrix \cite[\textsection 13A]{BurdeZieschang:1985}.
  % (As always, simplicity is partly a matter of taste.)
\end{example}

Section \ref{sub:Examples} displays some further examples 
to show that colouring polynomials distinguish knots 
for which other invariants fail:
\begin{itemize}
\item
  They distinguish the Kinoshita-Terasaka knot
  from the Conway knot and show that none of them is
  inversible nor reversible nor obversible.
\item
  They detect asymmetries of bretzel knots; 
  they distinguish, for example, $B(3,5,7)$
  from its inverse, reverse and obverse knot.
\item
  They distinguish the (inversible)
  knot $8_{17}$ from its reverse.
\end{itemize}

We also mention two natural questions that will not be pursued here:

\begin{question}
  Can knot colouring polynomials detect other geometric properties of knots?
  Applications to periodic knots and ribbon knots would be most interesting.
\end{question}

\begin{question}
  Do colouring polynomials distinguish all knots?
  Since the knot group system $(\pi_{K},m_K,l_K)$ charaterizes 
  the knot $K$ \cite[Cor.\,6.5]{Waldhausen:1968}, and knot groups 
  are residually finite \cite[Thm.\,3.3]{Thurston:1982}, 
  this question is not completely hopeless.
\end{question}

% \begin{question}
%   Given a knot $K$, how can we effectively find 
%   a suitable finite quotient $G$ of $\pi_{K}$?
%   The quotient group should be small enough 
%   to be manageable, but at the same time large enough 
%   to extract meaningful information about $K$.
% \end{question}

\subsection{Colouring polynomials are Yang-Baxter invariants}

Moving from empirical evidence to a more theoretical level, 
this article compares knot colouring polynomials with 
two other classes of knot invariants:
Yang-Baxter invariants, i.e.\  knot invariants 
obtained from Yang-Baxter representations of the braid group,
and quandle colouring state-sum invariants 
derived from quandle cohomology.
The result can be summarized as follows: 
\[
\Biggl\{ \begin{matrix}  \text{Yang-Baxter} \\ 
                         \text{invariants}         
\end{matrix} \Biggr\} 
\supset % \supsetneqq 
\Biggl\{ \begin{matrix}  \text{colouring} \\ 
                         \text{polynomials}        
\end{matrix} \Biggr\}
\supset % \supsetneqq 
\Biggl\{ \begin{matrix}  \text{quandle $2$-cocycle} \\ 
                         \text{state-sum invariants} 
\end{matrix} \Biggr\}
\supset % =
\Biggl\{ \begin{matrix}  \text{colouring} \\ 
                         \text{polynomials} \\
                         \text{with $\Lambda$ abelian} 
\end{matrix} \Biggr\}
\]

P.J.\,Freyd and D.N.\,Yetter \cite[Prop.\,4.2.5]{FreydYetter:1989} 
have shown that every colouring number $\F{x}{G}\colon \knots \to\Z$  
can be obtained from a certain Yang-Baxter operator $c$ over $\Z$.
We generalize this result to colouring polynomials:

\begin{theorem}[\textsection\ref{sub:YBclosed}] \label{Thm:YBclosed}
  Suppose that $G$ is a group with basepoint $x$
  such that the subgroup $\Lambda = C(x)\cap G'$ is abelian.
  Then the colouring polynomial $\P{x}{G}\colon \knots \to \Z\Lambda$ 
  is a Yang-Baxter invariant of closed knots: there exists 
  a Yang-Baxter operator $\tilde{c}$ over the ring $\Z\Lambda$, 
  such that the associated knot invariant 
  coincides with (a constant multiple) of $\P{x}{G}$.
\end{theorem}

In the general case, where $\Lambda$ is not necessarily abelian,
Section \ref{sub:YBlong} gives an analogous presentation of $\P{x}{G}$ 
as a Yang-Baxter invariant of long knots (also called $1$-tangles).

\begin{corollary}
  Since $\Lambda$ is abelian in all our examples of \textsection\ref{sub:Examples},
  it follows in particular that Yang-Baxter invariants 
  \emph{can} detect non-inversible and non-reversible knots.
\end{corollary}

\begin{remark}
  It follows from our construction that $\tilde{c}$ 
  is a deformation of $c$ over the ring $\Z\Lambda$.
  Conversely, the deformation ansatz leads to quandle 
  cohomology (see \textsection\ref{sub:Conclusion}).  
  Elaborating this approach, M.\,Gra\~na \cite{Grana:2002} 
  has shown that quandle $2$-cocycle state-sum invariants 
  are Yang-Baxter invariants.  The general theory 
  of Yang-Baxter deformations of $c_Q$ over 
  the power series ring $\K\mathopen{[\![} h \mathclose{]\!]}$ 
  has been developed in \cite{Eisermann:2005}.
\end{remark}  

\begin{remark}
  The celebrated Jones polynomial and, more generally, 
  all quantum invariants of knots, can be obtained 
  from Yang-Baxter operators that are formal 
  power series deformations of the trivial operator.  
  This implies that the coefficients in this expansion 
  are of finite type \cite[\textsection2.1]{BarNatan:1995}.
  Part of their success lies in the fact that these invariants distinguish 
  many knots, and in particular they easily distinguish mirror images. 
  It is still unknown, however, whether finite type invariants 
  can detect non-inversible or non-reversible knots.

  For colouring polynomials the construction is similar in that 
  $\P{x}{G}$ arises from a deformation of a certain operator $c$.
  There are, however, two crucial differences:
  \begin{itemize}
  \item 
    The initial operator $c$ models conjugation (and is not the trivial operator),
  \item 
    Its deformation $\tilde{c}$ is defined over $\Z\Lambda$ (and not over a power series ring).
  \end{itemize}
  As a consequence, the colouring polynomial $\P{x}{G}$ is not of finite type,
  nor are its coefficients, nor any other rational-valued invariant
  computed from it \cite{Eisermann:2000}.
\end{remark}

\subsection{Quandle invariants are specialized colouring polynomials}

A quandle, as introduced by D.\,Joyce \cite{Joyce:1982}, 
is a set $Q$ with a binary operation whose axioms model conjugation in a group,
or equivalently, the Reidemeister moves of knot diagrams.
Quandles have been intensively studied 
by different authors and under various names;
we review the relevant definitions in \textsection\ref{sec:Quandles}.
The Lifting Lemma proved in \textsection\ref{sub:LiftingLemma}
tells us how to pass from quandle to group colourings 
and back without any loss of information. 
On the level of knot invariants this implies the following result:

\begin{theorem}[\textsection\ref{sub:LiftingLemma}]
  Every quandle colouring number $\F{q}{Q}$ is the specialization 
  of some knot colouring polynomial $\P{x}{G}$.
\end{theorem}

Quandle cohomology was initially studied in order
to construct invariants in low-dimensional topology:
in \cite{CarterEtAl:1999,CarterEtAl:2003} it was shown 
how a $2$-cocycle $\lambda\in Z^2(Q,\Lambda)$
gives rise to a state-sum invariant of knots, 
$\CSSI{Q}{\lambda}\colon \knots \to \Z\Lambda$,
which refines the quandle colouring number $\TF{Q}$. % $\F{q}{Q}$.
We prove the following result:

\begin{theorem}[\textsection\ref{sub:SS2CP}] \label{Thm:SS2CP}
  Every quandle $2$-cocycle state-sum invariant of knots
  is the specialization of some knot colouring polynomial.
  More precisely, suppose that $Q$ is a connected quandle, 
  $\Lambda$ is an abelian group, and  $\lambda \in Z^2(Q,\Lambda)$ 
  is a $2$-cocycle with associated invariant 
  $\CSSI{Q}{\lambda}\colon \knots \to \Z{\Lambda}$.
  Then there exists a group $G$ with basepoint $x$ and 
  a $\Z$-linear map $\varphi\colon \Z{G} \to \Z{\Lambda}$ 
  such that $\CSSI{Q}{\lambda} = \varphi\P{x}{G} \cdot |Q|$.
\end{theorem}

This result provides a complete topological interpretation 
of quandle $2$-cocycle state-sum invariants in terms of 
the knot group and its peripheral system.
Conversely, we prove that state-sum invariants 
contain those colouring polynomials $\P{x}{G}$ for which 
the longitude group $\Lambda = C(x) \cap G'$ is abelian:

\begin{theorem}[\textsection\ref{sub:CP2SS}] \label{Thm:CP2SS}
  Suppose that $G$ is a colouring group with basepoint $x$
  such that the subgroup $\Lambda = C(x)\cap G'$ is abelian.
  Then the colouring polynomial $\P{x}{G}$ can be presented
  as a quandle $2$-cocycle state-sum invariant.
  More precisely, the quandle $Q = x^G$ admits 
  a $2$-cocycle $\lambda\in Z^2(Q,\Lambda)$
  such that $\CSSI{Q}{\lambda} = \P{x}{G} \cdot |Q|$.
\end{theorem}

% It should be noted, however, that knot colouring polynomials
% form a strictly larger class, because they also include cases 
% where $\Lambda$ is non-abelian (see Example \ref{exm:Alt9}).

\subsection{How this article is organized}

Section \ref{sec:KnotGroups} recalls the necessary facts
about the knot group and its peripheral system.
It then discusses connected sum and symmetry operations with 
respect to knot colouring polynomials and displays some applications.
The main purpose is to give some evidence as to 
the scope and the usefulness of these invariants.

Section \ref{sec:Quandles} examines quandle colourings
and %The Lifting Lemma (\textsection\ref{sub:LiftingLemma})
explains how to replace quandle colourings 
by group colourings without any loss of information.
The correspondence between quandle extensions and quandle cohomology
is then used to show how quandle $2$-cocycle state-sum invariants
can be seen as specializations of colouring polynomials.

Section \ref{sec:YangBaxter} relates colouring polynomials with 
Yang-Baxter invariants. After recalling the framework of linear
braid group representations, we show how colouring polynomials 
can be seen as Yang-Baxter deformations of colouring numbers.

\subsection{Acknowledgements}

The author would like to thank the anonymous referee
for his careful reading and numerous helpful comments.
The results of Section \ref{sec:KnotGroups} were part of 
the author's Ph.D.\ thesis \cite{Eisermann:Thesis}, 
which was financially supported by the Deutsche Forschungs\-gemeinschaft
through the Gra\-du\-ierten\-kolleg Mathe\-matik at the University of Bonn.
Sections \ref{sec:Quandles} and \ref{sec:YangBaxter} were elaborated while 
the author held a post-doc position at the \'Ecole Normale Sup\'erieure de Lyon, 
whose hospitality is gratefully acknowledged.

%%%%%%%%%%%%%%%%%%%%%%%%%%%%%%%%%%%%%%%%%%%%%%%%%%%%%%%%%%%%%%%%%%%%%%%%%%%%%

\section{Knot groups and colouring polynomials} \label{sec:KnotGroups}

This section collects some basic facts about the knot group 
and its peripheral system (\textsection\ref{sub:PeripheralSystem})
and their homomorphic images (\textsection\ref{sub:ColouringGroups}).
We explain how connected sum and symmetry operations affect 
the knot group system and how this translates to colouring polynomials 
(\textsection\ref{sub:KnotGroupSymmetries}).
We then display some examples showing that colouring polynomials 
are a useful tool in distinguishing knots where other invariants fail
(\textsection\ref{sub:Examples}).

\subsection{Peripheral system} \label{sub:PeripheralSystem}

We use fairly standard notation, which we recall
from \cite{Eisermann:2003} for convenience.
A \emph{knot} is a smooth embedding $k \colon \S^1\hookrightarrow\S^3$,
considered up to isotopy.  This is equivalent to considering
the oriented image $K = k(\S^1)$ in $\S^3$, again up to isotopy.
A \emph{framing} of $k$ is an embedding 
$f \colon \S^1 \times \D^2 \hookrightarrow \S^3$ such that $f|_{\S^1{\times}0} = k$.
As basepoint of the space $\S^3 \minus K$ we choose $p = f(1,1)$.
In the fundamental group $\pi_{K} := \pi_1(\S^3\minus K,p)$
we define the \emph{meridian} $m_K = [f|_{1{\times}\S^1}]$ 
and the \emph{longitude} $l_K = [f|_{\S^1{\times}1}]$.
Up to isotopy the framing is characterized 
by the linking numbers $\lk(K,m_K) \in \{\pm1\}$ 
and $\lk(K,l_K) \in \Z$, and all combinations are realized.
We will exclusively work with the \emph{standard framing}, 
characterized by the linking numbers $\lk(K,m_K) = +1$ and $\lk(K,l_K)=0$.

Up to isomorphism, the triple $(\pi_{K},m_K,l_K)$ is a knot invariant,
and even a complete invariant:  two knots $K$ and $K'$ are isotopic 
if and only if there is a group isomorphism $\phi \colon \pi_{K} \to \pi_{K'}$
with $\phi(m_K) = m_{K'}$ and $\phi(l_K) = l_{K'}$.
This is a special case of Waldhausen's theorem on sufficiently large 
$3$-manifolds; see \cite[Cor.\,6.5]{Waldhausen:1968} as well as 
\cite[\textsection{3C}]{BurdeZieschang:1985}.

Besides closed knots $k\colon \S^1 \hookrightarrow \S^3$ 
it will be useful to consider long knots (also called $1$-tangles), 
i.e.\  smooth embeddings $\ell\colon \R \hookrightarrow \R^3$ 
such that $\ell(t) = (t,0,0)$ for all parameters $t$ 
outside of some compact interval.  We refer to 
\cite{Eisermann:2003} for a detailed discussion
with respect to knot groups and quandles.

\begin{figure}[hbtp]
  \centering
  \includegraphics[width=0.9\linewidth]{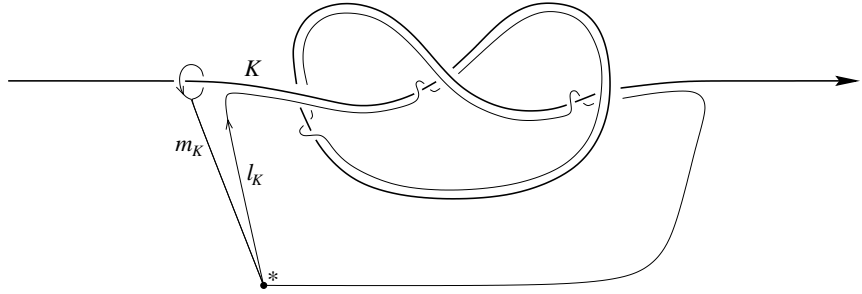}
  \caption{Meridian and longitude of a long knot}
  \label{fig:MeriLong}
\end{figure}

\subsection{Colouring groups} \label{sub:ColouringGroups}

Since knot groups are residually finite \cite[Thm.\,3.3]{Thurston:1982},
there are plenty of finite knot group representations.
But which groups do actually occur as homomorphic images of knot groups?
This question was raised by L.P.\,Neuwirth \cite{Neuwirth:1965},
and first solved by F.\,Gonzalez-Acu\~na \cite{Gonzalez:1975}:

\begin{theorem}[\cite{Gonzalez:1975,Johnson:1980}] 
  A pointed group $(G,x)$ is the homomorphic image of 
  some knot group $(\pi_{K},m_K)$ if and only if $G$ 
  is finitely generated and $G = \gen{x^G}$.
  \qed
\end{theorem}

The condition is necessary, because every knot group $\pi_{K}$
is finitely generated by conjugates of the meridian $m_K$.
(See the Wirtinger presentation, recalled in \textsection\ref{sub:WirtingerQuandles}.) 
For a proof of sufficiency we refer to the article of D.\,Johnson \cite{Johnson:1980},
who has found an elegant and ingeniously simple way to construct 
a knot $K$ together with an epimorphism $(\pi_{K},m_K) \to (G,x)$.
Here we restrict attention to \emph{finite} groups:

\begin{definition} \label{def:ColouringGroup}
  Let $G$ be a finite group and $x\in G$.
  The pair $(G,x)$ is called a \emph{colouring group}
  if the conjugacy class $x^G$ generates the whole group $G$.
  For example, every finite simple group $G$ is a colouring group
  with respect to any of its non-trivial elements $x \ne 1$.
\end{definition}

\begin{remark} \label{rem:ConnectedGroups}
  Given a finite group $G_0$ and $x \in G_0$,
  every homomorphism $(\pi_{K},m_K) \to (G_0,x)$
  maps to the subgroup $G_1 := \gen{x^{G_0}}$.
  If $G_1$ is strictly smaller than $G_0$,
  then we can replace $G_0$ by $G_1$.
  Continuing like this, we obtain a descending chain
  $G_0 \supset G_1 \supset G_2 \supset \cdots$,
  recursively defined by $G_{i+1} = \gen{x^{G_i}}$.
  Since $G_0$ is finite, this chain must stabilize, 
  and we end up with a colouring group $G_n = \gen{x^{G_n}}$.
  Hence, we can assume without loss of generality 
  that $(G,x)$ is a colouring group.
\end{remark}

Given $(G,x)$ let $\Lambda^*$ be the set of longitude images 
$\rho(l_K)$, where $\rho$ ranges over all knot group homomorphisms 
$\rho \colon (\pi_{K},m_K) \to (G,x)$ and all knots $K$.  
Then $\Lambda^*$ is a subgroup of $G$ \cite{JohnsonLivingston:1989}. 
Since meridian $m_K \in \pi_{K}$ and longitude $l_K \in \pi_{K}'$ commute,
$\Lambda^*$ is contained in the subgroup $\Lambda = C(x)\cap G'$,
which will play an important r\^ole in subsequent arguments.

D.\,Johnson and C.\,Livingston \cite{JohnsonLivingston:1989}
have worked out a complete characterization of 
the subgroup $\Lambda^*$ in terms of homological obstructions.
As an application, consider a colouring group $(G,x)$
that is perfect, i.e.\ $G'=G$, and has cyclic centralizer, 
say $C(x) = \gen{x}$.  Then \cite{JohnsonLivingston:1989} 
affirms that $\Lambda^* = \Lambda = C(x)$.  All of our examples 
in \textsection\ref{sub:Examples} are of this type.

\subsection{Knot and group symmetries} \label{sub:KnotGroupSymmetries}

The knot group $\pi_{K}$ is obviously independent of orientations.
In order to define the longitude, however, 
we have to specify the orientation of $K$,
and the definition of the meridian additionally
depends on the orientation of $\S^3$.
Changing these orientations defines
the following symmetry operations:

\begin{definition}
  Let $K\subset\S^3$ be an oriented knot.
  The same knot with the opposite orientation of $\S^3$
  is the \emph{mirror image} or the \emph{obverse} of $K$,
  denoted $K^\obv$. (We can represent this as $K^\obv = \sigma K$,
  where $\sigma\colon \S^3 \to \S^3$ is a reflection.)
  Reversing the orientation of the knot $K$
  yields the \emph{reverse} knot $K^\rev$.
  Inverting both orientations yields 
  the \emph{inverse} knot $K^\inv$.
\end{definition}

Please note that different authors use different terminology,
in particular reversion and inversion are occasionally interchanged.
Here we adopt the notation of J.H.\,Conway \cite{Conway:1969}.
% following the remark that the inverse knot $K^\inv$ represents 
% the inverse of $K$ in the knot concordance group \cite{FoxMilnor:1966}.
% Moreover, if $K$ is represented as the closure of a braid $\beta$
% (as explained in \textsection\ref{sub:YangBaxter}), then 
% the inverse braid $\beta^{-1}$ represents the inverse knot $K^\inv$.

% We derive the following easy observations;
% for details see \cite[\textsection 3C]{BurdeZieschang:1985}.

\begin{proposition}
  \newcommand{\GS}[1]{\check\pi(#1)}
  Let $K$ be an oriented knot with group system $\GS{K} = (\pi_{K},m_K,l_K)$.
  Obversion, reversion and inversion affect the group system as follows:
  \begin{alignat*}{3}
    & \text{obversion:} \quad && \GS{K^\obv} && = (\pi_{K},m_K^{-1},l_K) \\
    & \text{reversion:} \quad && \GS{K^\rev} && = (\pi_{K},m_K^{-1},l_K^{-1}) \\
    & \text{inversion:} \quad && \GS{K^\inv} && = (\pi_{K},m_K,l_K^{-1}) 
  \end{alignat*}
  The fundamental group of the connected sum $K \cs L$ 
  is the amalgamated product $\pi_{K} \ast \pi_{L}$ modulo $m_K=m_L$.
  Its meridian is $m_K$ and its longitude is the product $l_K l_L$.
  \qed
\end{proposition}

\begin{corollary} \label{cor:Multiplicative}
  Every colouring polynomial $\P{x}{G}\colon \knots \to \Z{G}$ is multiplicative, that is, 
  we have $\P{x}{G}(K\cs L) = \P{x}{G}(K) \cdot \P{x}{G}(L)$ for any two knots $K$ and $L$.
  \qed
\end{corollary}

In order to formulate the effect of inversion,
% on the colouring polynomial $\P{x}{G}$, 
let ${\;}^\inv\colon \Z{G} \to \Z{G}$ be the linear extension
of the inversion map $G \to G$, $g \mapsto g^{-1}$.

\begin{corollary} \label{cor:Inversion}
  Every colouring polynomial $\P{x}{G}\colon \knots \to \Z{G}$ is equivariant 
  under inversion, i.e.\  $\P{x}{G}(K^\inv) = \P{x}{G}(K)^\inv$ for every knot $K$.
  In particular, the colouring number $\F{x}{G}(K)$ is invariant under inversion of $K$.
  \qed
\end{corollary}

%\subsection{Group symmetries} \label{sub:GroupSymmetries}

Obversion and reversion of knots can similarly 
be translated into symmetries of colouring polynomials,
but to do so we need a specific automorphism of $G$:

\begin{definition} 
  % Let $G$ be a group with basepoint $x$.
  An automorphism ${\;}^\obv\colon G \to G$ with $x^\obv = x^{-1}$
  is called an \emph{obversion} of $(G,x)$.
  An anti-automorphism ${\;}^\rev\colon G \to G$ with $x^\rev = x$
  is called a \emph{reversion} of $(G,x)$.
\end{definition}

Obviously a group $(G,x)$ possesses a reversion 
if and only if it possesses an obversion.
They are in general not unique, because they can 
be composed with any automorphism $\alpha \in \Aut(G,x)$,
for example conjugation by an element in $C(x)$.

\begin{remark}
  The braid group $\Br{n}$, recalled in \textsection\ref{sub:YangBaxter}
  below, has a unique anti-automorphism ${\;}^\rev\colon \Br{n} \to \Br{n}$ 
  fixing the standard generators $\sigma_1,\dots,\sigma_{n-1}$.
  Analogously there exists a unique automorphism
  ${\;}^\obv\colon \Br{n} \to \Br{n}$ mapping each standard generator
  $\sigma_i$ to its inverse $\sigma_i^{-1}$.
  The exponent sum $\Br{n} \to \Z$ shows that 
  this cannot be an inner automorphism.
  % These two maps provide us with a reversion resp.\ obversion 
  % of the braid group $\Br{n}$ with basepoint $\sigma_1$.

  These symmetry operations on braids correspond 
  to the above symmetry operations on knots:
  if a knot $K$ is represented as the closure 
  of the braid $\beta$ (see \textsection\ref{sub:YangBaxter}),
  then the inverse braid $\beta^{-1}$ represents the inverse knot $K^\inv$,
  the reverse braid $\beta^\rev$ represents the reverse knot $K^\rev$,
  and the obverse braid $\beta^\obv$ represents the obverse knot $K^\obv$.
\end{remark}

Given an obversion and a reversion of $(G,x)$,
their linear extensions to the group ring $\Z{G}$
will also be denoted by ${\;}^\obv\colon \Z{G}\to\Z{G}$ 
and ${\;}^\rev\colon \Z{G}\to\Z{G}$, respectively.
We can now formulate the equivariance of 
the corresponding colouring polynomials:

\begin{corollary} \label{cor:CPEquivariance}
  Suppose that $(G,x)$ possesses an obversion 
  ${\;}^\obv$ and a reversion ${\;}^\rev$.
  Then the colouring polynomial $\P{x}{G}$ is equivariant
  with respect to obversion and reversion, 
  that is, we have $\P{x}{G}(K^\obv) = \P{x}{G}(K)^\obv$ 
  and  $\P{x}{G}(K^\rev) = \P{x}{G}(K)^\rev$ for every knot $K$.
  In this case the colouring numbers of 
  $K$, $K^\inv$, $K^\obv$, and $K^\rev$ are the same.
  \qed
\end{corollary}

% \begin{proof}
%   For a knot $K$ with group system $(\pi_{K},m_K,l_K)$,
%   the obverse knot $K^\obv$ has group system $(\pi_{K},m_K^{-1},l_K)$.
%   The automorphism ${\;}^\obv\colon (G,x) \to (G,x^{-1})$ establishes a bijection
%   between representations $\rho\colon (\pi_{K},m_K) \to (G,x)$ 
%   and representations $\rho^\obv\colon (\pi_{K},m_K^{-1}) \to (G,x)$,
%   hence $\P{x}{G}(K^\obv) = \P{x}{G}(K)^\obv$.
%   Reversion is the composition of obversion and inversion.
%   Equivariance then follows from $\P{x}{G}(K^\obv) = \P{x}{G}(K)^\obv$
%   and $\P{x}{G}(K^\inv) = \P{x}{G}(K)^\inv$.
% \end{proof}

\begin{example}
  Every element $x$ in the symmetric group $\Sym{n}$ is conjugated 
  to its inverse $x^{-1}$, because both have the same cycle structure. 
  Any such conjugation defines an obversion $(\Sym{n},x) \to (\Sym{n},x^{-1})$.
  This argument also applies to alternating groups: given $x \in \Alt{n}$ 
  we know that $x$ is conjugated to $x^{-1}$ in $\Sym{n}$.
  Since $\Alt{n}$ is normal in $\Sym{n}$, this conjugation 
  restricts to an obversion $(\Alt{n},x) \to (\Alt{n},x^{-1})$.
  This need not be an inner automorphism.
\end{example}

On the other hand, some groups do not permit any obversion at all:
% The simplest example is the following:

\begin{example}
  Let $\GF$ be a finite field and let 
  $G = \GF \rtimes \GF^\times$ be its affine group.
  We have $\Aut(G) = \Inn(G) \rtimes \Gal(\GF)$,
  where $\Gal(\GF)$ is the Galois group
  of $\GF$ over its prime field $\GF[p]$.
  % (The proof is left as an exercise.)
  If $\GF = \GF[p]$, then every automorphism of $G$ is inner and thus 
  induces the identity on the abelian quotient $\GF^\times$.
  If $p\ge5$, we can choose an element $x=(a,b) \in G$ % with $b\ne\pm1$.
  whose projection to $\GF^\times$ satisfies $b \ne b^{-1}$.
  Hence there is no automorphism of $G$ that maps $x$ to $x^{-1}$.
  Indeed, searching all groups of small order with GAP \cite{GAP},
  we find that the smallest group having this property 
  is $\GF[5] \rtimes \GF[5]^\times$ of order $20$.
\end{example}

% \begin{Notes}

For the sake of completeness we expound the following elementary result:

\begin{proposition}
  The affine group $G = \GF \rtimes \GF^\times$
  satisfies $\Aut(G) = \Inn(G) \rtimes \Gal(\GF)$.
\end{proposition}

\begin{proof}
  The product in $G$ is given by $(a,b)(c,d) = (a+bc,bd)$,
  and so $\Gal(\GF)$ can be seen as a subgroup of $\Aut(G)$,
  where $\phi \in \Gal(\GF)$ acts as $(a,b) \mapsto (\phi(a),\phi(b))$.
  Since $\Inn(G)$ is a normal subgroup of $\Aut(G)$
  with $\Inn(G) \cap \Gal(\GF) = \{\id_G\}$, we see that $\Aut(G)$ 
  contains the semi-direct product $\Inn(G) \rtimes \Gal(\GF)$.
  
  It remains to show that every $\alpha \in \Aut(G)$
  belongs to $\Inn(G) \rtimes \Gal(\GF)$.
  This is trivially true for $\GF = \GF[2]$,
  so we will assume that $\GF$ has more than two elements.
  It is then easily verified that $G' = \GF \times \{1\}$.
  Let $\zeta$ be a generator of the multiplicative group $\GF^\times$.
  We have $\alpha(1,1) = (u,1)$ with $u \in \GF^\times$, and 
  $\alpha(0,\zeta) = (v,\xi)$ with $v \in \GF$, $\xi \in \GF^\times$, $\xi \ne 1$.
  Conjugating by $w = ( v (1-\xi)^{-1}, u )$,
  we obtain $(u,1)^w = (1,1)$ and $(v,\xi)^w = (0,\xi)$.
  In the sequel we can thus assume $u = 1$ and $v = 0$.
  This implies $\alpha(0,b) = (0,\phi(b))$ 
  with $\phi \colon \GF^\times \to \GF^\times$, 
  $\zeta^n \mapsto \xi^n$ for all $n \in \Z$.
  Extending this by $\phi(0)=0$ we obtain a bijection $\phi \colon \GF \to \GF$
  satisfying $\phi(ab) = \phi(a)\phi(b)$ for all $a,b \in \GF$.
  Moreover, we find $\alpha(a,1) = (\phi(a),1)$:
  this is clear for $a = 0$, and for $a \ne 0$
  we have $(a,1) = (0,a)(1,1)$ and thus 
  $\alpha(a,1) = (0,\phi(a)) (1,1) = (\phi(a),1)$.
  This proves that $\phi(a+b) = \phi(a) + \phi(b)$ 
  for all $a,b \in \GF$, whence $\phi \in \Gal(\GF)$.  
  We conclude that $\alpha(a,b) = (\phi(a),\phi(b))$, as claimed.
\end{proof}

% \end{Notes}

\subsection{Examples and applications} \label{sub:Examples}

The preceding discussion indicates 
that symmetries of the group $(G,x)$
affect the colouring polynomial $\P{x}{G}(K)$
just as well as symmetries of the knot $K$.
We point out several examples:

\begin{example} \label{exm:PSLObv}
  Let $p$ be a prime and let $G=\PSL{2}{\GF[p]}$ be equipped with basepoint
  $z = \left[\begin{smallmatrix} 1 & 1 \\ 0 & 1 \end{smallmatrix}\right]$
  of order $p$. Inversion, obversion, and reversion are realized by
  \[
  \left[\begin{smallmatrix} a & b \\ c & d \end{smallmatrix}\right]^\inv = 
  % \left[\begin{smallmatrix} a & b \\ c & d \end{smallmatrix}\right]^{-1} = 
  \left[\begin{smallmatrix} d & -b \\ -c & a \end{smallmatrix}\right], \qquad
  \left[\begin{smallmatrix} a & b \\ c & d \end{smallmatrix}\right]^\obv = 
  \left[\begin{smallmatrix} a & -b \\ -c & d \end{smallmatrix}\right], \qquad
  \left[\begin{smallmatrix} a & b \\ c & d \end{smallmatrix}\right]^\rev = 
  \left[\begin{smallmatrix} d & b \\ c & a \end{smallmatrix}\right].
  \]
  We have $C(z) = \gen{z}$. For $p=2$ and $p=3$ one finds 
  that the longitude group $\Lambda = C(z)\cap G'$ is trivial.
  For $p \ge 5$ the group $G$ is perfect (even simple), hence $\Lambda = \gen{z}$.
  We conclude that the colouring polynomial $\P{z}{G}$ 
  is insensitive to reversion: we have $\P{z}{G}(K) \in \Z\gen{z}$
  and reversion fixes $z$ and therefore all elements in $\Z\gen{z}$.
\end{example}

\begin{example} \label{exm:AltObv}
  Consider an alternating group $G = \Alt{n}$ with $n \ge 3$,
  and a cycle $x = (123\dots l)$ of maximal length, 
  that is, $l=n$ for $n$ odd and $l=n-1$ for $n$ even.
  As we have pointed out above, a suitable conjugation 
  in $\Sym{n}$ produces an obversion $(G,x) \to (G,x^{-1})$.
  We have $C(x) = \gen{x}$.  For $n = 3$ and $n = 4$ one finds that 
  the longitude group $\Lambda = C(x) \cap G'$ is trivial. For $n\ge5$ the group $G$ 
  is perfect (even simple), hence the longitude group is $\Lambda = \gen{x}$.
  Again we conclude that the colouring polynomial $\P{x}{G}$ is insensitive to reversion.

  We observe that for $n=3,4,7,8,11,12,\dots$ an obversion 
  of $(G,x)$ cannot be realized by an inner automorphism: 
  consider for example $G = \Alt{11}$ and $x = (\mathtt{abcdefghijk})$:
  in $\Sym{11}$ the centralizer is $C(x) = \gen{x}$ and consequently every permutation 
  $\sigma \in \Sym{11}$ with $x^\sigma = x^{-1}$ is of the form 
  $\sigma = x^k (\mathtt{ak})(\mathtt{bj})(\mathtt{ci})(\mathtt{dh})(\mathtt{eg})$
  and thus odd.  The same argument shows that for $n=5,6,9,10,\dots$
  an obversion of $(G,x)$ can be realized by an inner automorphism.
\end{example}

\begin{example} \label{exm:MathieuObv}
  As a more exotic example, let us finally consider the Mathieu group $M_{11}$,
  i.e.\ the unique simple group of order $7920 = 2^4{\cdot}3^2{\cdot}5{\cdot}11$,
  and the smallest of the sporadic simple groups \cite{ATLAS}.
  It can be presented as a subgroup of $\Alt{11}$, 
  for example as 
  \[
  G = \gen{x,y} \quad\text{with}\quad 
  x = (\mathtt{abcdefghijk}),\; y=(\mathtt{abcejikdghf}).
  \]
  This presentation has been obtained from GAP \cite{GAP}
  and can easily be verified with any group-theory software
  by checking that $G$ is simple of order $7920$.
  % \begin{alltt}
  %   gap> G:= Group([(1,2,3,4,5,6,7,8,9,10,11),(1,2,3,5,10,9,11,4,7,8,6)]);; 
  %   gap> Size(G); IsSimple(G);
  %   7920
  %   true
  % \end{alltt}
  % According to \cite{ATLAS} there is only one simple group of order $7920$.
  The Mathieu group $M_{11}$ is particularly interesting for us,
  because it does \emph{not} allow an obversion. To see this it suffices 
  to know that its group of outer automorphisms is trivial \cite{ATLAS},
  in other words, every automorphism of $M_{11}$ is realized by conjugation.
  In $M_{11}$ the element $x$ is not conjugated to its inverse --- 
  this is not even possible in $\Alt{11}$ according to the preceding example.  
  Hence there is no automorphism of $M_{11}$ that maps $x$ to $x^{-1}$.
\end{example}

\begin{figure}[hbtp]
  \centering
  \includegraphics[height=30mm]{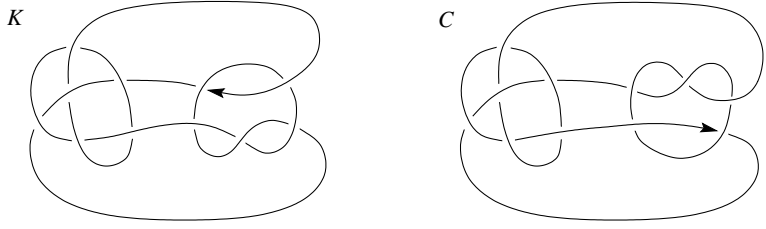}
  \caption{The Kinoshita-Terasaka knot and the Conway knot}
  \label{fig:ConwayKnot}
\end{figure}

Applied to colouring polynomials, this means that there is a priori 
no restriction on the invariants of a knot and its mirror image.
As a concrete example we consider the Kinoshita-Terasaka knot $K$
and the Conway knot $C$ displayed in Figure \ref{fig:ConwayKnot}.

Both knots have trivial Alexander polynomial.
They differ only by rotation of a $2$-tangle,
in other words they are mutants 
in the sense of Conway \cite{Conway:1969}.
Therefore neither the Jones, \textsc{Homflypt} nor Kauffman polynomial
can distinguish between $K$ and $C$, see \cite{Lickorish:1997}. 
With the help of a suitable colouring polynomial
the distinction is straightforward:

\begin{example} \label{exm:ConwaySL}
  R.\,Riley \cite{Riley:1971} has studied knot group homomorphisms
  to the simple group $G = \PSL{2}{\GF[7]}$ of order $168$.
  % of order $168 = 2^3\cdot3\cdot7$. % \cong \PSL{3}{\GF[2]}$ 
  Let $z$ be an element of order $7$, say 
  $z = \left[\begin{smallmatrix} 1 & 1 \\ 0 & 1 \end{smallmatrix}\right]$.
  % $z = \pm\left(\begin{smallmatrix} 1 & 1 \\ 0 & 1 \end{smallmatrix}\right)$.
  Then the associated colouring polynomials are
  \begin{alignat*}{3}
    % & \P{z}{G}(K)     & \;=\; & \P{z}{G}(C)      & \;=\; & 1 + 7z^{-1} + 7z^{-2} , \\
    & \P{z}{G}(K)      & \;=\; & \P{z}{G}(C)      & \;=\; & 1 + 7z^5 + 7z^6 , \\
    & \P{z}{G}(K^\inv) & \;=\; & \P{z}{G}(C^\inv) & \;=\; & 1 + 7z + 7z^2 .
  \end{alignat*}
  This shows that both knots are chiral.
  By a more detailed analysis of their coverings, 
  Riley could even show that $K$ and $C$ are distinct.
\end{example}

\begin{example}
  To distinguish $K$ and $C$ we give a simple and direct argument
  using colouring polynomials. For every element $x\in\PSL{2}{\GF[7]}$ of order $3$, 
  say $x = \left[\begin{smallmatrix} 0 & 1 \\ -1 & 1 \end{smallmatrix}\right]$,
  % $x = \pm\left(\begin{smallmatrix} 0 & 1 \\ -1 & 1 \end{smallmatrix}\right)$,
  the associated colouring polynomial distinguishes $K$ and $C$:
  \begin{alignat*}{4}
    & \P{x}{G}(K)      & \;=\; & 1 + 6x \hspace{2cm} &
    & \P{x}{G}(C)      & \;=\; & 1 + 12x \\
    % & \P{x}{G}(K^\inv) & \;=\; & 1 + 6x^{-1} \hspace{2cm} &
    % & \P{x}{G}(C^\inv) & \;=\; & 1 + 12x^{-1}
    & \P{x}{G}(K^\inv) & \;=\; & 1 + 6x^2 \hspace{2cm} &
    & \P{x}{G}(C^\inv) & \;=\; & 1 + 12x^2 
  \end{alignat*}
  Both invariants, $\P{z}{G}$ and $\P{x}{G}$,
  show chirality but are insensitive to reversion.
\end{example}

These and the following colouring polynomials were calculated with 
the help of an early prototype of the computer program \textsl{KnotGRep},
an ongoing programming project to efficiently construct 
the set of knot group homomorphisms to a finite group. 
Even though general-purpose software may be less comfortable,
our results can also be obtained from the Wirtinger 
presentation (\textsection\ref{sub:WirtingerQuandles})
using GAP \cite{GAP} or similar group-theoretic software.

\begin{example} \label{exm:ConwayAlt}
  The alternating group $G=\Alt{7}$ with basepoint $x=(1234567)$ yields
  \begin{alignat*}{4}
    & \P{x}{G}(K)      & \;=\; & 1 + 7x^2 + 28x^5 + 28x^6 \hspace{1cm} &
    & \P{x}{G}(C)      & \;=\; & 1 + 7x^2 + 7x^3 +21x^5 + 14x^6 \\
    & \P{x}{G}(K^\inv) & \;=\; & 1 + 28x + 28x^2 + 7x^5 \hspace{1cm} &
    & \P{x}{G}(C^\inv) & \;=\; & 1 + 14x + 21x^2 + 7x^4 + 7x^5 
  \end{alignat*}
  Again this invariant distinguishes $K$ et $C$ 
  and shows their chirality, but is insensitive to reversion,
  as explained in Example \ref{exm:AltObv} above.
\end{example}

\begin{example} \label{exm:ConwayMathieu}
  More precise information can be obtained using the Mathieu group $M_{11}$,
  presented as the permutation group $(G,x)$ in Example \ref{exm:MathieuObv} above.
  For the Kinoshita-Terasaka knot $K$ and the Conway knot $C$ one finds:
  \begin{alignat*}{4}
    & \P{x}{G}(K)      & \;=\; & 1 + 11x^3 + 11x^7 \hspace{2cm} &
    & \P{x}{G}(C)      & \;=\; & 1 + 11x^3 + 11x^7 \\
    & \P{x}{G}(K^\inv) & \;=\; & 1 + 11x^4 + 11x^8 \hspace{2cm} &
    & \P{x}{G}(C^\inv) & \;=\; & 1 + 11x^4 + 11x^8 \\
    & \P{x}{G}(K^\obv) & \;=\; & 1 + 11x^4 + 22x^8 \hspace{2cm} &
    & \P{x}{G}(C^\obv) & \;=\; & 1 + 11x^4 + 11x^6 + 11x^8 \\
    & \P{x}{G}(K^\rev) & \;=\; & 1 + 22x^3 + 11x^7 \hspace{2cm} &
    & \P{x}{G}(C^\rev) & \;=\; & 1 + 11x^3 + 11x^5 + 11x^7 
  \end{alignat*}
  Consequently all eight knots are distinct;
  $K$ and $C$ are neither inversible nor obversible nor reversible.
  (This example was inspired by G.\,Kuperberg \cite{Kuperberg:1996},
  who used the colouring number $\F{x}{G}$ to distinguish 
  the knot $C$ from its reverse $C^\rev$.)
\end{example}

\begin{figure}[hbtp]
  \centering
  \includegraphics[height=30mm]{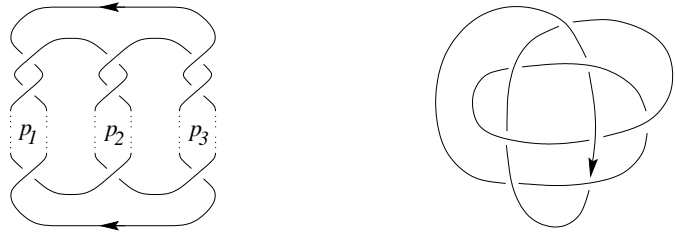}
  %\sidebyside{\caption{the bretzel knot $B(p_1,p_2,p_3)$}}{\caption{the knot $8_{17}$}}
  \caption{(a) the bretzel knot $B(p_1,p_2,p_3)$, (b) the knot $8_{17}$}
  \label{fig:Bretzel}
\end{figure}

Usually it is very difficult to detect non-reversibility of knots.
Most invariants fail to do so, including the usual knot polynomials.
In view of the simplicity of our approach, the success
of knot colouring polynomials is remarkable.
We give two further examples:

\begin{example} \label{exm:BretzelMathieu}
  The family of bretzel knots $B(p_1,p_2,p_3)$,
  parametrized by odd integers $p_1,p_2,p_3$,
  is depicted in Figure \ref{fig:Bretzel}a.
  According to the classification of bretzel knots 
  (see \cite{BurdeZieschang:1985}, \textsection 12),
  the bretzel knot $B=B(3,5,7)$ is neither 
  reversible nor obversible nor inversible.
  For the Mathieu group $G=M_{11}$ with basepoint $x$
  as in Example \ref{exm:ConwayMathieu} we obtain:
  \begin{alignat*}{4}
    & \P{x}{G}(B)      & \;=\; & 1 + 11x \hspace{2cm} &
    & \P{x}{G}(B^\obv) & \;=\; & 1 + 11x^7 \\
    & \P{x}{G}(B^\inv) & \;=\; & 1 + 11x^{10} \hspace{2cm} &
    & \P{x}{G}(B^\rev) & \;=\; & 1 + 11x^4
  \end{alignat*}
  Again the colouring polynomial shows that the knot $B$
  possesses none of the three symmetries.
  Historically, bretzel knots were 
  the first examples of non-reversible knots.
  Their non-reversibility was first proven 
  by H.F.\,Trotter \cite{Trotter:1963} in 1963
  by representing the knot group system on a suitable 
  triangle group acting on the hyperbolic plane.
\end{example}

\begin{example} \label{exm:817Mathieu}
  Figure \ref{fig:Bretzel}b shows the knot $8_{17}$, 
  which is the smallest non-reversible knot.
  It is a $3$-bridge knot but not a bretzel knot,
  and there is no general classification theorem available.
  To analyze this example we choose once more
  the Mathieu group $M_{11}$ with basepoint $x$ as above.
  The knot $8_{17}$ then has colouring polynomial $1+11x^5+11x^6$
  whereas the reverse knot has trivial colouring polynomial $1$.
  (Here even the colouring number $\F{x}{G}$ 
  suffices to prove that this knot is non-reversible.) 
  We remark that $8_{17}$ is inversible and that this symmetry 
  is reflected in the symmetry of its colouring polynomials.
\end{example}

The colouring polynomial $\P{x}{G}(K)$ is, 
by definition, an element in the group ring $\Z{G}$,
and it actually lies in the much smaller ring $\Z\Lambda$.
The following symmetry consideration 
further narrows down the possible values.
It is included here to explain one of the observations 
that come to light in the previous examples, 
but it will not be used in the sequel.

\begin{proposition} \label{prop:PrimeMultiple}
  Let $(G,x)$ be a colouring group.
  If conjugation by $x$ has order $p^k$ for some prime $p$, then 
  the colouring polynomial satisfies $\P{x}{G}(K) \equiv 1 \pmod{p}$.
\end{proposition}

\begin{proof}
  The cyclic subgroup $\gen{x}$ acts on the set 
  $\Hom(\pi_{K},m_K \,;\, G,x)$ by conjugation.
  The only fixed point is the trivial representation 
  $(\pi(K),m_K) \to (\Z,1) \to (G,x)$.  This can be most easily
  seen by interpreting group homomorphisms $\rho \colon (\pi_{K},m_K) \to (G,x)$ 
  as colourings $f \colon (D,0) \to (G,x)$ of a knot diagram $D$, 
  see \textsection\ref{sub:WirtingerQuandles} below.
  If $f^x = f$ then all colours of $f$ commute with $x$: 
  following the diagram from the first to the last arc 
  we see by induction that all colours are in fact equal to $x$.
  Since there is only one component, we conclude 
  that $f$ is the trivial colouring, corresponding 
  to the trivial representation. 

  Every non-trivial representation $\rho$ appears 
  in an orbit of length $p^\ell$ for some $\ell\ge1$.
  Since $\rho(l_K)$ commutes with $x$, all representations 
  in such an orbit have the same longitude image in $G$.
  The sum $\P{x}{G}(K)$ thus begins with $1$ for the trivial representation,
  and all other summands can be grouped to multiples of $p$.
\end{proof}

\section{Quandle invariants are specialized colouring polynomials} \label{sec:Quandles}

The Wirtinger presentation allows us to interpret 
knot group homomorphisms as colourings of knot diagrams.
Since such colourings involve only conjugation, 
they are most naturally treated in the category of quandles, 
as introduced by D.\,Joyce \cite{Joyce:1982}.
We recall the basic definitions concerning quandles 
and quandle colourings in \textsection\ref{sub:WirtingerQuandles},
and explain in \textsection\ref{sub:LiftingLemma} 
how to pass from quandles to groups and back
without any loss of information.

Quandle cohomology was studied in \cite{CarterEtAl:1999,CarterEtAl:2003}, 
where it was shown how a $2$-cocycle gives rise to 
a state-sum invariant of knots in $\S^3$.
We recall this construction in \textsection\ref{sub:CP2SS}
and show that every colouring polynomial $\P{x}{G}$
can be presented as a quandle $2$-cocycle state-sum invariant,
provided that the subgroup $\Lambda=C(x)\cap G'$ is abelian
(Theorem \ref{thm:CP2SS}).

In order to prove the converse, we employ the 
cohomological classification of central quandle extensions 
established in \cite{Eisermann:2003,CarterEtAl:2003b},
recalled in \textsection\ref{sub:QuandleCoverings} below.
This allows us to prove in \textsection\ref{sub:SS2CP}
that every quandle $2$-cocycle state-sum invariant is the specialization 
of a suitable knot colouring polynomial (Theorem \ref{thm:SS2CP}).
% This result provides a complete topological interpretation 
% of quandle $2$-cocycle state-sum invariants in terms of 
% the knot group and its peripheral system.

\subsection{Wirtinger presentation, quandles, and colourings} \label{sub:WirtingerQuandles}

\begin{figure}[hbtp]
  \centering
  \includegraphics[height=20mm]{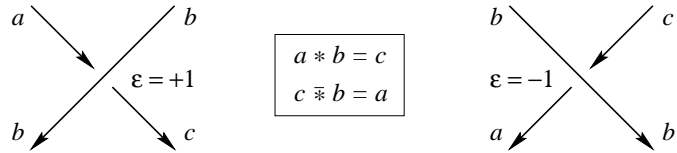}
  \caption{Wirtinger rules for colouring a knot diagram}
  \label{fig:ColouredCrossing}
\end{figure}

Our exposition follows \cite{Eisermann:2003},
to which we refer for further details.
We consider a long knot diagram as in Figure \ref{fig:MeriLong}
and number the arcs consecutively from $0$ to $n$.
At the end of arc number $i-1$, we undercross 
arc number $\kappa{i}=\kappa(i)$ and continue on arc number $i$.
We denote by $\varepsilon{i}=\varepsilon(i)$ the sign of 
this crossing, as depicted in Figure \ref{fig:ColouredCrossing}.
The maps $\kappa\colon \{1,\dots,n\}\to\{0,\dots,n\}$ 
and $\varepsilon\colon \{1,\dots,n\}\to\{\pm1\}$
are the \emph{Wirtinger code} of the diagram.

\begin{theorem} \label{thm:Wirtinger}
  Suppose that a knot $L$ is represented by a long knot diagram 
  with Wirtinger code $(\kappa,\varepsilon)$ as above.
  Then the knot group allows the presentation
  \[
  \pi_{L} = \langle x_0,x_1,\dots,x_n \,|\, r_1,\dots,r_n \rangle
  \quad\text{with relation $r_i$ being}\quad 
  x_i = x_{\kappa{i}}^{-\varepsilon{i}} \, x_{i-1} \, x_{\kappa{i}}^{\varepsilon{i}} .
  \]
  As peripheral system we can choose $m_L = x_0$ and 
  $l_L = \prod_{i=1}^{i=n} \; x_{i-1}^{-\varepsilon{i}} \, x_{\kappa{i}}^{\varepsilon{i}}$.
  \qed
\end{theorem}

For a proof see Crowell-Fox \cite[\textsection{VI.3}]{CrowellFox:1963}
or Burde-Zieschang \cite[\textsection{3B}]{BurdeZieschang:1985}.
The Wirtinger presentation works just as well for a closed knot diagram.
Since arcs $0$ and $n$ are then identified, this amounts to adding 
the (redundant) relation $x_0=x_n$ to the above presentation. 
The group is, of course, the same.

The Wirtinger presentation allows us to interpret 
knot group homomorphisms $\pi_{L} \to G$ as colourings.
More precisely, a \emph{$G$-colouring} of the diagram $D$ 
is a map $f \colon \{0,\dots,n\} \to G$ such that 
$f(i) = f(\kappa{i})^{-\varepsilon{i}} f(i-1) f(\kappa{i})^{\varepsilon{i}}$.
In other words, at each coloured crossing as in 
Figure \ref{fig:ColouredCrossing} the colours $a$ and $c$ 
are conjugated via $a^b=c$.  Such a colouring is denoted by $f \colon D \to G$.
We denote by $\Col(D;G)$ the set of colourings of $D$ with colours in $G$.
For a long knot diagram $D$, we denote by $\Col(D,0;G,x)$
the subset of colourings that colour arc number $0$ with colour $x$.
The Wirtinger presentation establishes natural bijections
$\Hom(\pi_{K};G) \cong \Col(D;G)$ and $\Hom(\pi_{K},m_K;G,x) \cong \Col(D,0;G,x)$.

\begin{figure}[hbtp]
  \centering
  \includegraphics[width=0.8\textwidth]{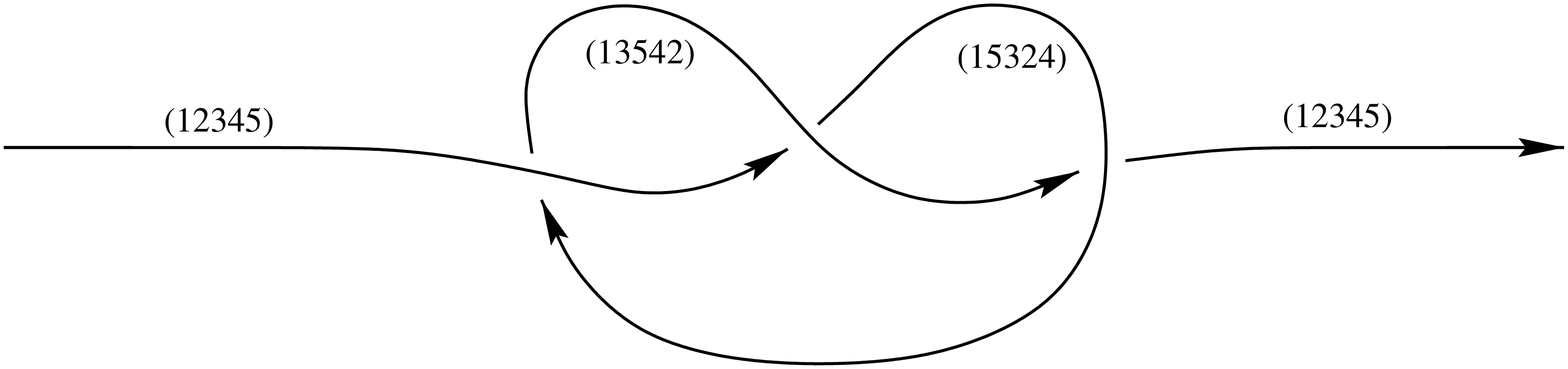} %height=30mm
  \caption{$\Alt{5}$-colouring of the left-handed trefoil knot}
  \label{fig:TrefoilAlt5}
\end{figure}

\begin{example}
  Figure \ref{fig:TrefoilAlt5} shows a colouring of the left-handed trefoil knot
  (represented as a long knot) with elements in the alternating group $\Alt{5}$.
  Note that all definitions readily extend to closed knot diagrams.
\end{example}

%\subsection{Quandles and colourings} \label{sub:Quandles}

The Wirtinger presentation of $\pi_{K}$ involves 
only conjugation but not the group multiplication itself.
The underlying algebraic structure can be described as follows:

\begin{definition}
  A \emph{quandle} is a set $Q$ with two binary operations 
  $\ast,\tsa\colon Q\times Q\to Q$ satisfying the following axioms
  for all $a,b,c\in Q$:
  \begin{enumerate} \addtolength{\itemindent}{\parindent}
  \item[(Q1)]
    $a\ast a = a$ \hfill (idempotency)
  \item[(Q2)]
    $(a \ast b) \tsa b = (a \tsa b) \ast b = a$ \hfill (right invertibility)
  \item[(Q3)]
    $(a\ast b)\ast c = (a\ast c)\ast (b\ast c)$ \hfill (self-distributivity)
  \end{enumerate}
\end{definition}

The name ``quandle'' was introduced by D.\,Joyce \cite{Joyce:1982}.
The same notion was studied by S.V.\,Matveev \cite{Matveev:1982}
under the name ``distributive groupoid'', and by L.H.\,Kauffman 
\cite{Kauffman:2001} who called it ``crystal''.
Quandle axioms (Q2) and (Q3) are equivalent to saying that for every $b \in Q$ 
the right translation $\varrho_b \colon a \mapsto a \ast b$ is an automorphism of $Q$.
Such structures were called ``automorphic sets'' by E.\,Brieskorn \cite{Brieskorn:1988}.
The somewhat shorter term \emph{rack} was preferred by R.\,Fenn and 
C.P.\,Rourke \cite{FennRourke:1992}.  The notion has been generalized
to ``crossed $G$-sets'' by P.J.\,Freyd and D.N.\,Yetter \cite{FreydYetter:1989}.

\begin{definition}
  As before, let $D$ be a long knot diagram, 
  its arcs being numbered by $0,\dots,n$.
  A \emph{$Q$-colouring}, denoted $f\colon D\to Q$, 
  is a map $f\colon \{0,\dots,n\}\to Q$ such that 
  at each crossing as in Figure \ref{fig:ColouredCrossing} 
  the three colours $a,b,c$ satisfy the relation $a \ast b=c$.
  We denote by $\Col(D;Q)$ the set of $Q$-colourings,
  and by $\Col(D,0;Q,q)$ the subset of colourings satisfying $f(0)=q$.
\end{definition}

\begin{proposition}[Joyce \cite{Joyce:1982}]
  The quandle axioms ensure that each Reidemeister move $D \leftrightharpoons D'$ 
  induces bijections $\Col(D;Q) \leftrightharpoons \Col(D';Q)$ and
  $\Col(D,0;Q,q) \leftrightharpoons \Col(D',0;Q,q)$, respectively.
  % For details see Joyce \cite{Joyce:1982}, \textsection{15}. 
  In particular, if $Q$ is finite, then the colouring numbers $\TF{Q}(D) = |\Col(D;Q)|$ 
  and $\F{q}{Q}(D) = |\Col(D,0;Q,q)|$ are knot invariants.
  \qed
\end{proposition}

\subsection{From quandle colourings to group colourings and back} \label{sub:LiftingLemma}

In many respects quandles are close to groups. 
For colourings we will now explain how to pass 
from quandles to groups and back without any loss of information.
% We will now have a closer look at this relationship 
% in order to explain how to pass from quandles to groups and back.

\begin{definition}
  A \emph{quandle homomorphism} is a map $\phi\colon Q\to Q'$ 
  that satisfies $\phi(a\ast b) = \phi(a) \ast \phi(b)$,
  and hence $\phi(a\tsa b) = \phi(a) \tsa \phi(b)$, for all $a,b\in Q$.
\end{definition}

\begin{definition}
  The automorphism group $\Aut(Q)$ consists of 
  all bijective homomorphisms $\phi\colon Q\to Q$.
  We adopt the convention that automorphisms of $Q$ act on the right,
  written $a^\phi$, which means that their composition $\phi\psi$ 
  is defined by $a^{(\phi\psi)} = (a^\phi)^\psi$ for all $a\in Q$.
\end{definition}
  
\begin{definition}
  The group $\Inn(Q) = \gen{ \varrho_b \mid b \in Q}$ 
  of \emph{inner automorphisms} is the subgroup of $\Aut(Q)$ 
  generated by all right  translations $\varrho_b \colon a \mapsto a \ast b$.
  The quandle $Q$ is called \emph{connected} 
  if the action of $\Inn(Q)$ on $Q$ is transitive.
\end{definition}

% As for colouring groups, we can restrict attention
% to connected quandles, see Remark \ref{rem:ConnectedGroups}.

In view of the map $\varrho \colon Q \to \Inn(Q)$, $b \mapsto \varrho_b$,
we also write $a^b = a \ast b$ for the operation in a quandle.
Conversely, it will sometimes be convenient to write 
$a \ast b = b^{-1} a b$ for the conjugation in a group.
In neither case will there be any danger of confusion.

\begin{definition}
  A \emph{representation} of a quandle $Q$ 
  on a group $G$ is a map $\phi\colon Q\to G$ such that 
  $\phi(a\ast b) = \phi(a) \ast \phi(b)$ for all $a,b\in Q$.
  In other words, the following diagram commutes:
  \[
  \begin{CD}
    Q \times Q @>{\phi \times \phi}>> G \times G \\
    @V{*}VV @VV{\conj}V \\
    Q @>{\phi}>> G 
  \end{CD}
  \]
\end{definition}

For example, the natural map $\varrho\colon Q \to \Aut(Q)$ 
satisfies $\varrho(a \ast b) = \varrho(a) \ast \varrho(b)$.
We call $\varrho$ the \emph{inner} representation of $Q$.
Moreover it satisfies $\varrho(a^g) = \varrho(a)^g$ 
for all $a \in Q$ and $g \in \Aut(Q)$.
This is the prototype of an augmentation:

\begin{definition}
  Let $\phi \colon Q \to G$ be a representation and let 
  $\alpha \colon Q \times G \to Q$, $(a,g) \mapsto a^g$,
  be a group action.  We call the pair $(\phi,\alpha)$
  an \emph{augmentation} if $a \ast b = a^{\phi(b)}$ and
  $\phi(a^g) = \phi(a)^g$ for all $a,b \in Q$ and $g \in G$.
  In other words, the following diagram commutes:
  \begin{equation} \label{eq:Augmentation}
  \begin{CD}
    Q \times Q @>{\id \times \phi}>> Q \times G @>{\phi \times \id}>> G \times G \\
    @V{*}VV @VV{\alpha}V @VV{\conj}V \\
    Q @>{\id}>> Q @>{\phi}>> G 
  \end{CD}
  \end{equation}
\end{definition}

\begin{remark}
  We will usually reinterpret the group action $\alpha$ 
  as a group homomorphism $\bar\alpha \colon G \to \Aut(Q)$,
  and denote the augmentation by $Q \lto[\phi] G \lto[\bar\alpha] \Aut(Q)$.
  If $G$ is generated by the image $\phi(Q)$, then $\phi$ is equivariant and 
  the action of $G$ on $Q$ is uniquely determined by the representation $\phi$.
  In this case we simply say that $\phi \colon Q \to G$ is an augmentation.
  For example, every quandle $Q$ comes equipped with 
  the inner augmentation $\varrho \colon Q \to \Inn(Q)$.
\end{remark}

Suppose that $Q$ is a quandle and $\phi\colon Q \to G$ 
is a representation on some group $G$.  Obviously every
quandle colouring $\tilde{f}\colon D \to Q$ maps to a
group colouring $f = \phi\tilde{f} \colon D \to G$.
If $\phi$ is an augmentation, then this process can be reversed,
and we can replace quandle colourings by group colourings 
without any loss of information:

\begin{lemma} \label{lem:LiftingLemma}
  Let $\phi \colon (Q,q) \to (G,x)$ be an augmentation 
  % of the quandle $Q$ on the group $G$ such that $\phi(q)=x$.
  of the quandle $Q$ with basepoint $q \in Q$ 
  on the group $G$ with basepoint $x \in G$.
  If $D$ is a long knot diagram, then every group colouring $f \colon (D,0) \to (G,x)$ 
  can be lifted to a unique quandle colouring $\tilde{f} \colon (D,0) \to (Q,q)$.
  %% satisfying $f = \phi\tilde{f}$.  
  In other words, $\phi$ induces a bijection 
  \[
  \phi_* \colon \Col(D,0;Q,q) \lto[\sim] \Col(D,0;G,x) ,
  \quad \tilde{f} \mapsto f = \phi\tilde{f} .
  \]
  Moreover, let $\rho \colon (\pi_{K},m_K) \to (G,x)$ be the knot group representation
  associated with $f$. Then the lifted colouring $\tilde{f}$ begins 
  with $\tilde{f}(0)=q$ and ends with $\tilde{f}(n) = q^{\rho(l_K)}$.
\end{lemma}

\begin{proof}
  Every representation $\phi \colon (Q,q) \to (G,x)$ induces a map $\phi_*$ 
  sending each quandle colouring $\tilde{f} \colon (D,0) \to (Q,q)$ to 
  the associated group colouring $\phi\tilde{f} \colon (D,0) \to (G,x)$.  
  In general $\phi_*$ is neither injective nor surjective, 
  lest $\phi$ is an augmentation. In order to define the inverse map 
  $\psi_*\colon \Col(D,0;G,x) \to \Col(D,0;Q,q)$, we use the action 
  $\alpha \colon Q \times G \to Q$, which we temporarily denote 
  by $(a,g) \mapsto a \bullet g$ for better readability.
  
  The crucial ingredient in the proof is 
  the commutativity of the diagram \eqref{eq:Augmentation}.
  Let us first show how the condition $a \ast b = a \bullet \phi(b)$
  ensures injectivity of $\phi_*$.  Let $D$ be a long knot 
  diagram with Wirtinger code $(\kappa,\varepsilon)$.
  Assume that $\tilde{f},\hat{f} \colon (D,0) \to (Q,q)$ 
  are colourings with $\phi\tilde{f} = \phi\hat{f}$.  
  By hypothesis we have $\tilde{f}(0) = \hat{f}(0) = q$.  
  By induction suppose that $\tilde{f}(i-1) = \hat{f}(i-1)$
  for some $i \ge 1$.  In the case of a positive crossing 
  ($\varepsilon i = +1$) we then obtain
  \begin{align*}
  \tilde{f}(i) 
  & = \tilde{f}(i-1) \ast \tilde{f}(\kappa i)
  = \tilde{f}(i-1) \bullet \phi\tilde{f}(\kappa i)
  \\
  & = \hat{f}(i-1) \bullet \phi\hat{f}(\kappa i)
  = \hat{f}(i-1) \ast \hat{f}(\kappa i)
  = \hat{f}(i) .
  \end{align*}
  The case of a negative crossing ($\varepsilon i = -1$) is analogous.

  We now show how the equivariance condition $\phi(a \bullet g) = \phi(a) \ast g$
  of diagram \eqref{eq:Augmentation} ensures surjectivity.
  For every colouring $f \colon (D,0) \to (G,x)$, denoted by $i \mapsto x_i$, 
  the colours $x_0,\dots,x_n$ satisfy $x_i = x_{i-1} \ast x_{\kappa{i}}^{\varepsilon{i}}$.
  We define partial longitudes $\ell_0,\dots,\ell_n$ by setting
  $
  \ell_i := \prod_{j=1}^{i} x_{j-1}^{-\varepsilon{j}} x_{\kappa{j}}^{\varepsilon{j}}.
  $
  In particular we have $x_0 = x_n = x$ and $x_i = x_0 \ast \ell_i$ for all $i=0,\dots,n$.
  By definition, $\ell_n = \rho(l_K)$ is the (total) longitude of the colouring $f$.
  We define $\tilde{f} \colon (D,0) \to (Q,q)$ by assigning 
  the colour $q_i = q \bullet \ell_i$ to arc number $i=0,\dots,n$.
  By hypothesis, $\phi \colon Q \to G$ is an equivariant map, whence 
  \begin{equation}
    \label{eq:ProjectionProperty}
    \phi(q_i) = \phi(q \bullet \ell_i) = \phi(q) \ast \ell_i = x \ast \ell_i = x_i.
  \end{equation}
  At each positive crossing we find the 
  following identity, using axiom (Q1):
  \begin{equation}
    \label{eq:PositiveCrossing}
    q_{i-1} \ast q_{\kappa{i}} 
    = (q_{i-1} \tsa q_{i-1}) \ast q_{\kappa{i}} 
    % = q_{i-1} \bullet (x_{i-1}^{-1} x_{\kappa{i}}^{\phantom{-1}}) 
    = (((q \bullet \ell_{i-1}^{\phantom{-1}}) \bullet x_{i-1}^{-1}) \bullet x_{\kappa{i}}^{\phantom{-1}}
    % = q \bullet (\ell_{i-1}^{\phantom{-1}} x_{i-1}^{-1} x_{\kappa{i}}^{\phantom{-1}})
    = q \bullet \ell_i = q_i.
  \end{equation}
  Analogously at each negative crossing:
  \begin{equation}
    \label{eq:NegativeCrossing}
    q_{i-1} \tsa q_{\kappa{i}} 
    = (q_{i-1} \ast q_{i-1}) \tsa q_{\kappa{i}}  
    % = q_{i-1} \bullet (x_{i-1}^{\phantom{-1}} x_{\kappa{i}}^{-1}) 
    = (((q \bullet \ell_{i-1}^{\phantom{-1}}) \bullet x_{i-1}^{\phantom{-1}}) \bullet x_{\kappa{i}}^{-1})
    % = q \bullet (\ell_{i-1}^{\phantom{-1}} x_{i-1}^{\phantom{-1}} x_{\kappa{i}}^{-1})
    = q \bullet \ell_i = q_i.
  \end{equation}
  
  % This shows that $\tilde{f} \colon (D,0) \to (Q,q)$ is a colouring.  
  We can thus define $\psi_* \colon \Col(D,0;G,x) \to \Col(D,0;Q,q)$ 
  by $f \mapsto \tilde{f}$. % by sending $f$ to its lifted colouring $\tilde{f}$. 
  Equation \eqref{eq:ProjectionProperty} shows that $\phi_*\psi_* = \id$, 
  and Equations \eqref{eq:PositiveCrossing} and \eqref{eq:NegativeCrossing} 
  imply that $\psi_*\phi_* = \id$.
\end{proof}

\begin{remark}
  Obviously, the condition $a \ast b = a^{\phi(b)}$ 
  cannot be dropped because it connects the quandle 
  operation $\ast$  with the group action $\alpha$.
  Likewise, the equivariance condition $\phi(a^g) = \phi(a)^g$ 
  cannot be dropped:  as an extreme counter-example, consider 
  a trivial quandle $Q = \{q\}$ and an arbitrary group $(G,x)$.
  We have a unique representation $\phi \colon (Q,q) \to (G,x)$
  and a unique group action $\alpha \colon Q \times G \to Q$.
  The map $\phi$ is equivariant if and only if $x \in Z(G)$.
  In general $\phi_*$ cannot be a bijection, because 
  the only $(Q,q)$-colouring is the trivial one, 
  while there may be non-trivial $(G,x)$-colourings.
\end{remark}

% The Lifting Lemma can be interpreted as saying that the longitude 
% acts as the \emph{monodromy} of a colouring along a knot.
% This point of view will become even more evident
% when we consider quandle coverings in \textsection\ref{sec:QuandleExtensions}.

The Lifting Lemma has the following analogue for closed knots:

\begin{lemma} \label{lem:ClosedLifting}
  Let $\phi \colon (Q,q) \to (G,x)$ be an augmentation 
  of the quandle $Q$ on the group $G$.
  If $D$ is a closed knot diagram, 
  then $\phi$ induces a bijection between 
  $\Col(D,0;Q,q)$ and those homomorphisms 
  $\rho \colon (\pi_{K},m_K) \to (G,x)$ 
  satisfying $q^{\rho(l_K)} = q$.
  \qed
\end{lemma}

As an immediate consequence we obtain the following result:

\begin{theorem}
  Every quandle colouring number $\F{q}{Q}$ is the specialization 
  of some knot colouring polynomial $\P{x}{G}$.
\end{theorem}

\begin{proof}
  We consider an augmentation $\phi \colon (Q,q) \to (G,x)$ 
  with $G = \gen{\phi(Q)}$, for example the inner augmentation 
  on $\phi \colon Q \to G = \Inn(Q)$ with basepoint $x = \phi(q)$. 

  For long knots, Lemma \ref{lem:LiftingLemma} implies $\F{q}{Q} = \F{x}{G}$.
  Hence $\F{q}{Q} = \varepsilon\P{x}{G}$, where $\varepsilon\colon \Z{G}\to\Z$
  is the augmentation map of the group ring, with $\varepsilon(g) = 1$ for all $g \in G$.

  For closed knots we define the linear map 
  $\varepsilon \colon \Z{G}\to\Z$ by setting $\varepsilon(g) = 1$ 
  if $q^g = q$, and $\varepsilon(g) = 0$ if $q^g \ne q$.
  Then Lemma \ref{lem:ClosedLifting} implies that
  $\F{q}{Q} = \varepsilon \P{x}{G}$.
\end{proof}

This argument will be generalized in \textsection\ref{sub:SS2CP},
where we show that every quandle $2$-cocycle state-sum invariant
is the specialization of some colouring polynomial.

\subsection{Quandle coverings, extensions, and cohomology} \label{sub:QuandleCoverings}

Following \cite{Eisermann:2003} we recall how quandle colourings 
can be used to encode longitudinal information.
To this end we consider a long knot diagram with 
meridians $x_0,\dots,x_n$ and partial longitudes $l_0,\dots,l_n$
as defined in the above proof of the lifting lemma.
In particular we have $x_0 = x_n = m_K$ and $x_i = x_0 \ast l_i$ 
with $l_0=1$ and $l_n=l_K$. 
If we colour each arc not only with its meridian $x_i$ 
but with the pair $(x_i,l_i)$, then at each crossing we find that
\[
x_i = x_{i-1} \ast x_{\kappa{i}}^{\varepsilon{i}} \quad\text{and}\quad
l_i = l_{i-1} x_{i-1}^{-\varepsilon{i}} x_{\kappa{i}}^{\varepsilon{i}} .
\]

This crossing relation can be encoded in a quandle as follows.

\begin{lemma}[\cite{Eisermann:2003}] \label{lem:CoveringQuandle}
  Let $G$ be a group that is generated by a conjugacy class $Q=x^G$.
  Then $Q$ is a connected quandle with respect to conjugation 
  $a \ast b = b^{-1} a b$ and its inverse $a \tsa b = b a b^{-1}$.
  Let $G'$ be the commutator subgroup and define
  \[
  \tilde{Q} = \tilde{Q}(G,x) := \{\; (a,g)\in G\times G' \;|\; a=x^g \;\}.
  \]
  The set $\tilde{Q}$ becomes a connected quandle when equipped with the operations 
  \[
  (a,g) \ast (b,h) = ( a \ast b, ga^{-1}b ) \quad\text{and}\quad
  (a,g) \tsa (b,h) = ( a \tsa b, gab^{-1} ).
  \]
  Moreover, the projection $p\colon \tilde{Q}\to Q$ given by 
  $p(a,g)=a$ is a surjective quandle homomorphism.
  It becomes an equivariant map when we let $G'$ act on $Q$ 
  by conjugation and on $\tilde{Q}$ by $(a,g)^b = (a^b,gb)$.
  In both cases $G'$ acts transitively and 
  as a group of inner automorphisms.
  \qed
\end{lemma}

The construction of the quandle $\tilde{Q}(G,x)$ 
has been tailor-made to capture longitude information.
Considered purely algebraically, it is a covering
in the following sense:

\begin{definition} \label{def:Covering}
  A surjective quandle homomorphism $p\colon \tilde{Q}\to Q$ 
  is called a \emph{covering} if $p(\tilde{x}) = p(\tilde{y})$
  implies $\tilde{a}\ast\tilde{x} = \tilde{a}\ast\tilde{y}$ 
  for all $\tilde{a},\tilde{x},\tilde{y}\in\tilde{Q}$.
  In other words, the inner representation 
  $\tilde{Q}\to\Inn(\tilde{Q})$ factors through $p$. 
  This property allows us to define an action of $Q$ on $\tilde{Q}$ by setting 
  $\tilde{a}\ast x := \tilde{a}\ast\tilde{x}$ with $\tilde{x}\in p^{-1}(x)$. 
  % Note that every augmentation $\phi\colon \tilde{Q}\to G$ defines a covering 
  % $\tilde{Q}\to Q$ when restricted to its image $Q=\phi(\tilde{Q})$.
\end{definition}

In the construction of Lemma \ref{lem:CoveringQuandle},
the projection $p\colon \tilde{Q}\to Q$ is a covering map.
Moreover, covering transformations are given by the left action of 
$\Lambda = C(x)\cap G'$ defined by $\lambda\cdot(a,g) = (a,\lambda g)$.
This action satisfies the following axioms:

\begin{enumerate} \addtolength{\itemindent}{\parindent}
\item[(E1)]
  $(\lambda\tilde{x}) \ast \tilde{y} = \lambda (\tilde{x}\ast\tilde{y})$
  and $\tilde{x} \ast (\lambda\tilde{y}) = \tilde{x}\ast\tilde{y}$
  for all $\tilde{x},\tilde{y}\in\tilde{Q}$ and $\lambda\in\Lambda$.
\item[(E2)]
  $\Lambda$ acts freely and transitively on each fibre $p^{-1}(x)$.
\end{enumerate}

Axiom (E1) is equivalent to saying that $\Lambda$ acts by automorphisms 
and the left action of $\Lambda$ commutes with the right action of $\Inn(\tilde{Q})$.
We denote such an action by $\Lambda \curvearrowright \tilde{Q}$.
In this situation the quotient $Q := \Lambda\backslash\tilde{Q}$ 
carries a unique quandle structure that turns the projection 
$p \colon \tilde{Q} \to Q$ into a quandle covering.

\begin{definition} \label{def:CentrExt}
  An \emph{extension} $E\colon \Lambda \curvearrowright \tilde{Q} \to Q$
  consists of a surjective quandle homomorphism $\tilde{Q}\to Q$ 
  and a group action $\Lambda \curvearrowright \tilde{Q}$ 
  satisfying axioms (E1) and (E2).
  We call $E$ a \emph{central extension} if $\Lambda$ is abelian.
\end{definition}

Quandle extensions are an analogue of group extensions,
and central quandle extensions come as close 
as possible to imitating central group extensions.
Analogous to the case of groups, %central group extensions,
central quandle extensions are classified by 
the second cohomology group $H^2(Q,\Lambda)$,
see \cite{Eisermann:2003,CarterEtAl:2003b}. More precisely:

\begin{theorem}[\cite{Eisermann:2003}] \label{thm:CentrExtClass}
  Let $Q$ be a quandle, let $\Lambda$ be an abelian group,
  and let $\Ext(Q,\Lambda)$ be the set of equivalence classes 
  of central extensions of $Q$ by $\Lambda$.
  Given a central extension $E\colon \Lambda \curvearrowright \tilde{Q} \to Q$,
  each section $s\colon Q\to\tilde{Q}$ defines a $2$-cocycle $\lambda\colon Q \times Q \to \Lambda$.
  If $s'$ is another section, then the associated $2$-cocycle 
  $\lambda'$ differs from $\lambda$ by a $2$-coboundary.
  The map $E \mapsto [\lambda]$ so constructed induces
  a natural bijection $\Ext(Q,\Lambda) \cong H^2(Q,\Lambda)$.
  \qed
\end{theorem}

The relevant portion of the cochain complex 
$C^1 \lto[\delta^1] C^2 \lto[\delta^2] C^3$
is formed by $n$-cochains $\lambda\colon Q^n \to \Lambda$ satisfying 
$\lambda(a_1,\dots,a_n)=0$ whenever $a_i=a_{i+1}$ for some index $i$,
and the first two coboundary operators %are given by
$\delta^1(\mu)(a,b) = \mu(a) - \mu(a^b)$ and
$\delta^2(\lambda)(a,b,c) = \lambda(a,c) - \lambda(a,b) + \lambda(a^c,b^c) - \lambda(a^b,c)$.
For details, see \cite{CarterEtAl:1999,CarterEtAl:2003,Eisermann:2003}

\subsection{From colouring polynomials to state-sum invariants} \label{sub:CP2SS}

Let $D$ be a knot diagram and let $f$ be a colouring of $D$ with colours in $Q$.
Suppose that $\Lambda$ is an abelian group, written multiplicatively,
and that $\lambda \colon Q^2\to \Lambda$ is a $2$-cocycle.
For each coloured crossing $p$ as in Figure \ref{fig:ColouredCrossing}, 
we define its \emph{weight} by $\weight{\lambda}{p} := \lambda(a,b)^\varepsilon$.
The total weight of the colouring $f$ is the product 
$\weight{\lambda}{f} := \prod_p \weight{\lambda}{p}$ over all crossings $p$.
The \emph{state-sum} of the diagram $D$ is defined to be
$\CSSI{Q}{\lambda}(D) := \sum_f \; \weight{\lambda}{f}$,
where the sum in $\Z\Lambda$ is taken over all colourings $f\colon D \to Q$.
We recall the following results:

\begin{lemma}[\cite{CarterEtAl:1999,CarterEtAl:2003}]
  The state-sum $\CSSI{Q}{\lambda}$ is invariant under Reidemeister moves
  and thus defines a knot invariant $\CSSI{Q}{\lambda}\colon \knots \to \Z\Lambda$.
  \qed
\end{lemma}

\begin{lemma}[{\cite[Prop.\,4.5]{CarterEtAl:2003}}] \label{lem:GaugeInvariance}
  If the colouring $f\colon D\to Q$ is closed, that is $f(0) = f(n)$,
  then the weight $\weight{\lambda}{f}$ is invariant under addition of coboundaries.
  As a consequence, the state sum $\CSSI{Q}{\lambda}$ of a closed knot
  depends only on the cohomology class $[\lambda]$.
  \qed
\end{lemma}

\begin{lemma}[{cf.\ \cite[Lem.\,32]{Eisermann:2005}}] \label{lem:ConjugationInvariance}
  The diagonal action of $\Inn(Q)$ on $Q^n$ 
  induces the trivial action on $H^*(Q,\Lambda)$.
  As a consequence, for each closed colouring $f \colon D \to Q$ 
  and every inner automorphism $g\in\Inn(Q)$ we have 
  $\weight{\lambda}{f^g} = \weight{^{g\!}\lambda}{f} = \weight{\lambda}{f}$.
  \qed
\end{lemma}

This last result is well-known in group cohomology,
cf.\ Brown \cite[Prop.\,II.6.2]{Brown:1994}.
It seems to be folklore in quandle cohomology,
but I could not find a written account of it.
The necessary argument is provided by \cite[Lem.\,32]{Eisermann:2005}
in the more general setting of Yang-Baxter cohomology,
which immediately translates to Lemma \ref{lem:ConjugationInvariance}

% \begin{proof}
%   For $\lambda\in Z^n(Q,\Lambda)$ and $g \in \Inn(Q)$
%   we have to prove that $\lambda - {}^{g\!}\lambda \in B^n(Q,\Lambda)$.
%   It suffices to show this for $g = \varrho(a)$ with $a \in Q$,
%   because $\Inn(Q)$ is generated by $\varrho(Q)$.
%   In this case the cocycle condition for $\lambda$ 
%   implies $\lambda - {}^{g\!}\lambda = \delta\mu$ with 
%   $\mu(a_1,\dots,a_{n-1}) = (-1)^n \lambda(a_1,\dots,a_{n-1},a)$,
%   and we are done.
% \end{proof}

\begin{lemma}[{\cite[Lem.\,50]{Eisermann:2003}}] \label{lem:LiftingMonodromy} 
  Let $p \colon (\tilde{Q},\tilde{q}) \to (Q,q)$ be a central quandle extension.
  Given a long knot diagram $D$, every colouring $f \colon (D,0) \to (Q,q)$
  uniquely lifts to a colouring $\tilde{f} \colon (D,0) \to (\tilde{Q},\tilde{q})$
  such that $f = p \tilde{f}$.  
  If $f$ is closed then $\tilde{f}(n) = \weight{\lambda}{f} \cdot \tilde{q}$,
  where $[\lambda] \in H^2(Q,\Lambda)$ is the cohomology class
  associated with the extension $p$.
  \qed
\end{lemma}

These preliminaries being in place, we can now 
prove that every colouring polynomial $\P{x}{G}$ 
can be presented as a $2$-cocycle state-sum invariant,
provided that the subgroup $\Lambda=C(x)\cap G'$ is abelian.

\begin{theorem} \label{thm:CP2SS}
  Suppose that $G$ is a colouring group with basepoint $x$
  such that the subgroup $\Lambda = C(x)\cap G'$ is abelian.
  Then the colouring polynomial $\P{x}{G}$ can be presented 
  as a quandle $2$-cocycle state-sum invariant.
  More precisely, the quandle $Q = x^G$ admits
  a $2$-cocycle $\lambda\in Z^2(Q,\Lambda)$
  such that $\CSSI{Q}{\lambda} = \P{x}{G}\cdot|Q|$.
\end{theorem}

\begin{proof}
  Let $Q=x^G$ be the conjugacy class of $x$ in the group $G$,
  and let $\tilde{Q} = \tilde{Q}(G,x)$ be the covering quandle
  constructed in Lemma \ref{lem:CoveringQuandle}. Since $\Lambda$ is abelian,
  we obtain a central extension $\Lambda\curvearrowright\tilde{Q}\to Q$.
  Let $[\lambda]\in H^2(Q,\Lambda)$ be the associated cohomology class.
  As basepoints we choose $q=x$ in $Q$ and $\tilde{q} = (x,1)$ in $\tilde{Q}$.
  
  Let $D$ be a long diagram of some knot $K$, 
  let $f \colon (D,0) \to (Q,q)$ be a colouring, 
  let $\rho \colon (\pi_{K},m_K) \to (G,x)$ be 
  the corresponding knot group homomorphism, and let 
  $\tilde{f} \colon (D,0) \to (\tilde{Q},\tilde{q})$ be the lifting of $f$.
  On the one hand we have $\tilde{f}(n) = (x,\weight{\lambda}{f})$
  from Lemma \ref{lem:LiftingMonodromy}.
  On the other hand we have $\tilde{f}(n) = (x,\rho(l_K))$
  from the Wirtinger presentation.
  Thus $\rho(l_K) = \weight{\lambda}{f}$, and summing over all
  colourings $f \colon (D,0) \to (Q,q)$ yields $\P{x}{G}(K)$.
  
  To obtain the state-sum $\CSSI{Q}{\lambda}$ 
  we have to sum over all colourings $f \colon D\to Q$.  
  We have $\Col(D,Q) = \bigcup_{q'\in Q} \Col(D,0;Q,q')$.
  Since $Q$ is connected, for each $q'\in Q$ there exists $g\in G$
  such that $q^g = q'$.  Hence $f \mapsto f^g$ establishes 
  a bijection between $\Col(D,0;Q,q)$ and $\Col(D,0;Q,q')$.
  By Lemma \ref{lem:ConjugationInvariance} we have 
  $\weight{\lambda}{f} = \weight{\lambda}{f^g}$.
  Thus the state-sum over all colourings $f \colon (D,0) \to (Q,q')$ 
  again yields $\P{x}{G}$.  We conclude that
  $\CSSI{Q}{\lambda}(K) = \P{x}{G}(K) \cdot |Q|$.
\end{proof}

\subsection{From state-sum invariants to colouring polynomials} \label{sub:SS2CP}

Theorem \ref{thm:CP2SS} has the following converse,
which allows us to express quandle $2$-cocycle state-sum 
invariants by knot colouring polynomials.

\begin{theorem} \label{thm:SS2CP}
  Every quandle $2$-cocycle state-sum invariant of knots
  is the specialization of some knot colouring polynomial.
  More precisely, suppose that $Q$ is a connected quandle, 
  $\Lambda$ is an abelian group, and $\lambda \in Z^2(Q,\Lambda)$ 
  is a $2$-cocycle with associated invariant 
  $\CSSI{Q}{\lambda}\colon \knots\to\Z{\Lambda}$.
  Then there exists a group $G$ with basepoint $x$ and
  a linear map $\varphi\colon \Z{G}\to\Z{\Lambda}$ such that 
  the colouring polynomial $\P{x}{G}\colon \knots\to\Z{G}$
  satisfies $\CSSI{Q}{\lambda} = \varphi\P{x}{G} \cdot |Q|$.
\end{theorem}

\begin{proof}
  % Given a connected quandle $Q$ and a cocycle $\lambda\in Z^2(Q,\Lambda)$, 
  % we will produce a colouring group $(G,x)$ and a linear map 
  % $\varphi\colon \Z{G}\to\Z{\Lambda}$ as claimed.
  
  We first construct a suitable group $(G,x)$ together with 
  a linear map $\varphi\colon \Z{G}\to\Z{\Lambda}$.
  Let $\Lambda \curvearrowright \tilde{Q} \lto[p] Q$ be 
  the central extension associated with the $2$-cocycle $\lambda$,
  as explained in Theorem \ref{thm:CentrExtClass}.
  We put $G:= \Inn(\tilde{Q})$. The inner representation 
  $\tilde\varrho\colon \tilde{Q} \to G$ defines an augmented quandle 
  in the sense of \textsection\ref{sub:LiftingLemma}.
  We choose a basepoint $\tilde{q}\in\tilde{Q}$
  and set $x:=\tilde{\varrho}(\tilde{q})$. 
  % This defines our group $(G,x)$.

  We choose $q = p(\tilde{q})$ as basepoint of $Q$.
  Let $s\colon Q\to\tilde{Q}$ be a section that realizes the $2$-cocycle $\lambda$.
  Since $p$ is a covering, we obtain a representation 
  $\varrho\colon Q \to G$ by $\varrho = \tilde\varrho\compose s$.
  Conversely, we can define an action of $G$ on $Q$ 
  by setting $a^g = p( s(a)^g )$.  This turns 
  the representation $\varrho \colon Q \to G$ into an augmentation 
  and $p \colon \tilde{Q} \to Q$ into an equivariant map.
  Our notation being in place, we can now define the linear map 
  \[
  \varphi\colon \Z{G}\to\Z\Lambda \quad\text{by setting}\quad 
  \varphi(g) = \begin{cases} 0 & \text{if } q^g \ne q , \\
    \ell  & \text{if $q^g = q$ and $\ell \in \Lambda$ 
      such that $\tilde{q}^g = \ell\cdot\tilde{q}$.}
  \end{cases}
  \]
  
  It remains to prove that $\CSSI{Q}{\lambda} = \varphi\P{x}{G}\cdot|Q|$.
  Let $K$ be a knot represented by a long knot diagram $D$.
  The Lifting Lemma \ref{lem:ClosedLifting} grants us a bijection 
  between closed colourings $f\colon (D,0) \to (Q,q)$ and those homomorphisms 
  $\rho \colon (\pi_{K},m_K) \to (G,x)$ that satisfy $q^{\rho(l_K)} = q$. 
  Regarding the covering $\tilde{Q}$, we claim that
  $\tilde{q}^{\rho(l_K)}  = \weight{\lambda}{f} \cdot \tilde{q}$.
  To see this, let $\tilde{f} \colon (D,0) \to (\tilde{Q},\tilde{q})$ be the lifting of $f$.
  On the one hand we can apply the Lifting Lemma \ref{lem:ClosedLifting}
  to the augmentation $\tilde{Q}\to G$, 
  which yields $\tilde{f}(n) = \tilde{q}^{\rho(l_K)}$.
  On the other hand we can apply Lemma \ref{lem:LiftingMonodromy},
  which yields $\tilde{f}(n) = \weight{\lambda}{f} \cdot \tilde{q}$.
  
  The map $\varphi$ thus specializes 
  the knot colouring polynomial $\P{x}{G}(K)$ 
  to the state-sum $\sum_f \weight{\lambda}{f}$, at least
  if we restrict the summation to colourings $f \colon (D,0) \to (Q,q)$.
  Since $Q$ is connected, any other basepoint $q'$ yields 
  the same state-sum by Lemma \ref{lem:ConjugationInvariance}. 
  Summing over all $q'\in Q$, we thus obtain
  $\CSSI{Q}{\lambda} = \varphi\P{x}{G}\cdot|Q|$, as claimed.
\end{proof}

%%%%%%%%%%%%%%%%%%%%%%%%%%%%%%%%%%%%%%%%%%%%%%%%%%%%%%%%%%%%%%%%%%%%%%%%%%%%%

\section{Colouring polynomials are Yang-Baxter invariants} \label{sec:YangBaxter}

P.J.\,Freyd and D.N.\,Yetter \cite{FreydYetter:1989} have shown that
the colouring number $\F{x}{G}\colon \knots\to\Z$ is a Yang-Baxter invariant.
This means that $\F{x}{G}$ can be obtained as the trace of a linear 
braid group representation arising from a suitable Yang-Baxter operator $c$.

In this section we will show that the colouring polynomial 
$\P{x}{G}\colon \knots\to\Z\Lambda$ is also a Yang-Baxter invariant, 
obtained from a certain Yang-Baxter operator $\tilde{c}$ defined below.
It will follow from our construction that $\tilde{c}$ 
is a deformation of $c$ over $\Z\Lambda$.

\subsection{Braid group representations and Yang-Baxter invariants} \label{sub:YangBaxter}

The notion of Yang-Baxter invariants rests on two 
classical theorems: Artin's presentation of the braid groups 
and the Alexander-Markov theorem, which we will now recall.
Our exposition closely follows \cite{Eisermann:2005}
and is included here for convenience.

% Braid groups were introduced by E.\,Artin in 1925, 
% see \cite{Artin:1947,Birman:1974}.
% The following well-known theorem gives 
% a presentation by generators and relations.

\begin{theorem}[E.\,Artin \cite{Artin:1947}]
  The braid group on $n$ strands can be presented as
  \[
  \Br{n}  = \left\langle\; \sigma_1,\dots,\sigma_{n{-}1} \;\Big|\;
    \begin{matrix} 
      \hfill{} \sigma_i \sigma_j = \sigma_j \sigma_i  \hfill{} \quad \text{for } |i-j|\ge2 \\
      \sigma_i \sigma_j \sigma_i  = \sigma_j  \sigma_i \sigma_j \quad \text{for } |i-j|=1
    \end{matrix} \;\right\rangle,
  \]
  where the braid $\sigma_i$ performs a positive half-twist of the strands $i$ and $i+1$.
\end{theorem}

\begin{definition}
  Let $\K$ be a commutative ring and $V$ a $\K$-module.
  A \emph{Yang-Baxter operator} (or \emph{R-matrix})
  is an automorphism $c\colon V \tensor V \to V \tensor V$
  that satisfies the \emph{Yang-Baxter equation},
  also called \emph{braid relation}:
  \[
  (c \tensor\id_V)(\id_V\tensor c)(c \tensor\id_V) = 
  (\id_V\tensor c)(c \tensor\id_V)(\id_V\tensor c) 
  \qquad\text{in}\quad \Aut_\K(V^{\tensor3}).
  \]
  Here and in the sequel tensor products are 
  taken over $\K$ if no other ring is indicated.
\end{definition}

\begin{corollary}
  Given a Yang-Baxter operator $c$ and some integer $n\ge2$, we can define 
  automorphisms $c_i\colon V^{\tensor n} \to V^{\tensor n}$ by setting 
  \[
  c_i = \id_V^{\tensor(i{-}1)} \,\tensor\, c \,\tensor\, \id_V^{\tensor(n{-}i{-}1)}
  \quad\text{for } i=1,\dots,n-1.
  \]
  The Artin presentation implies that there exists, 
  for each $n$, a unique braid group representation
  $\rho_c^n\colon \Br{n} \to \Aut_\K(V^{\tensor n})$
  defined by $\rho_c^n(\sigma_i)=c_i$.
  \qed
\end{corollary}

\begin{figure}[hbtp]
  \centering
  \includegraphics[width=100mm]{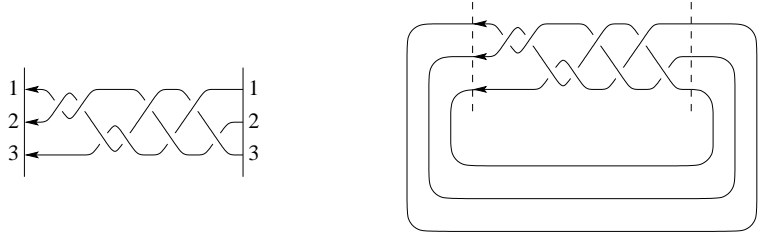}
  \caption{A braid $\beta$ and its closure $[\beta]$}
  \label{fig:BraidClosure}
\end{figure}

% In order to arrive at a consistent notation for braids 
% and their representations, we adopt the following conventions.
We orient braids from right to left as in Figure \ref{fig:BraidClosure}.
Braid groups will act on the left, so that composition of braids 
corresponds to the usual composition of maps.
% For example, the braid in Figure \ref{fig:BraidClosure} reads 
% $\beta = \sigma_1^{-2}\sigma_2^{+2}\sigma_1^{-1}\sigma_2^{+1}\sigma_1^{-1}\sigma_2^{+1}$.
% It is represented by the operator $\rho_c^3(\beta) = 
% c_1^{-2} c_2^{+2} c_1^{-1} c_2^{+1} c_1^{-1} c_2^{+1}$ acting on $V^{\tensor 3}$.
The passage from braids to links is granted by the closure map 
$[\,]\colon \bigcup_n \Br{n}\to\links$ defined as follows:
for each braid $\beta$ we define its \emph{closure} $[\beta]$ to be 
the link in $\S^3$ obtained by identifying opposite endpoints,
as indicated in Figure \ref{fig:BraidClosure}.

\begin{theorem}[Alexander-Markov, see \cite{Birman:1974}]
  Every link can be represented as the closure of some braid.
  Two braids represent the same link if and only if 
  one can be transformed into the other by a finite sequence 
  of the following \emph{Markov moves}:
  \begin{enumerate} \addtolength{\itemindent}{\parindent}
  \item[(M1)] 
    Pass from $\beta\in\Br{n}$ to 
    $\beta\sigma_n^{\pm1}\in\Br{n+1}$, or vice versa.
    \hfill \emph{(Stabilization)}
  \item[(M2)] 
    Pass from $\beta\in\Br{n}$ to % its conjugate
    $\alpha^{-1}\beta\alpha$ with $\alpha\in\Br{n}$.
    \hfill \emph{(Conjugation)}
  \end{enumerate}
\end{theorem}

Constructing a link invariant $F\colon \links\to\K$ is thus equivalent 
to constructing a map $F\colon \bigcup_n \Br{n}\to\K$ that is 
invariant under Markov moves.  The most natural approach is 
to consider traces of linear braid group representations:
invariance under conjugation is automatic, so we only have 
to require invariance under stabilization:

\begin{definition}
  Suppose that $V$ is a free $\K$-module with finite basis.
  Let $c \colon V \tensor V \to V \tensor V$ be a Yang-Baxter 
  operator.  An automorphism $m \colon V \to V$ is called 
  \emph{Markov operator} for $c$ if it satisfies 
  \begin{enumerate} \addtolength{\itemindent}{\parindent}
  \item[(m1)] 
    \makebox[35mm][l]{the trace condition}
    $ \tr_2 ( \; (m \tensor m) \compose c^{\pm1} \; ) = m$ \qquad and
  \item[(m2)] 
    \makebox[35mm][l]{commutativity}
    $c \compose (m \tensor m) = (m \tensor m) \compose c$.
  \end{enumerate}
\end{definition}

Here the partial trace $\tr_2 \colon \End(V \tensor V) \to \End(V)$ 
is defined as follows.  Let $(v_1,\dots,v_n)$ be a basis of $V$ over $\K$.
Every $f \in \End(V \tensor V)$ uniquely corresponds to a matrix $\smash{f_{ij}^{k\ell}}$ 
such that $f( v_i \tensor v_j ) = \sum_{k,\ell} \smash{f_{ij}^{k\ell}} \; v_k \tensor v_\ell$.
We can then define $g = \tr_2(f) \in \End(V)$, $g(v_i) = \sum_k \smash{g_i^k} \, v_k$, 
by the matrix $\smash{g_i^k} = \sum_j \smash{f_{ij}^{kj}}$.
(See Kassel \cite[\textsection II.3]{Kassel:1995}.)

\begin{corollary}
  Given a Yang-Baxter operator $c$ with Markov operator $m$,
  we define a family of maps $F_n\colon \Br{n}\to\K$ by 
  $F_n(\beta) = \tr( m^{\tensor n} \compose \rho_c^n(\beta) )$.
  Then the induced map $F\colon \bigcup_n \Br{n}\to\K$ 
  is invariant under both Markov moves and 
  thus defines a link invariant $F\colon \links\to\K$.
  \qed
\end{corollary}

The proof of this corollary is straight-forward:
the trace condition (m1) implies invariance under stabilization (M1),
and commutativity (m2) implies invariance under conjugation (M2).
Much more intricate is the question how to actually 
\emph{find} such a Yang-Baxter-Markov operator $(c,m)$.
Attempts to construct solutions in a systematic way have led 
to the theory of quantum groups \cite{Drinfeld:1987}.
For details we refer to the concise introduction 
\cite{KasselRossoTuraev:1997} or the textbook \cite{Kassel:1995}.

\begin{remark}
  For some Yang-Baxter operators $c$ there 
  does not exist any Markov operator $m$ at all.
  If it exists, $m$ is in general not the identity, 
  as in the case of the Jones polynomial 
  or other quantum invariants.
  The Yang-Baxter operators derived 
  from knot diagram colourings below
  are very special in that they allow 
  the Markov operator $m=\id$, which is 
  equivalent to saying that $\tr_2(c^{\pm1}) = \id$.
\end{remark}

\subsection{Colouring polynomials of long knots} \label{sub:YBlong}

Before we consider colouring polynomials, 
let us first recall how colouring numbers 
can be obtained from a suitable Yang-Baxter operator.
The following result is due to Freyd and Yetter, 
see \cite{FreydYetter:1989}, Prop.\ 4.2.5 and 
the remark following its proof.
% The following result is due to Freyd and Yetter for crossed $G$-sets,
% but its reformulation for quandles is straightforward:

\begin{theorem}[\cite{FreydYetter:1989}] \label{thm:CN2YB}
  Let $Q$ be a quandle and let $\K{Q}$ be the free $\K$-module with basis $Q$.
  The quandle structure of $Q$ can be linearly extended to a Yang-Baxter operator 
  \[
  c_Q\colon \K{Q}\tensor\K{Q} \to \K{Q}\tensor\K{Q} \quad\text{with}\quad
  a\tensor b \mapsto b \tensor (a\ast b) \quad\text{for all}\quad a,b\in Q.
  \]
  Axiom (Q2) ensures that $c_Q$ is an automorphism,
  while Axiom (Q3) implies the Yang-Baxter equation.
  If $Q$ is finite, then (Q1) ensures that $\tr_2(c_Q^{\pm1})=\id$.
  In this case the corresponding Yang-Baxter invariant 
  $\TF{Q} = \tr\compose\rho_Q$ coincides with the number 
  of $Q$-colourings (defined in \textsection\ref{sub:WirtingerQuandles})
  followed by the ring homomorphism $\Z \to \K$.
  \qed
\end{theorem}

As an example consider a finite group $G$ with basepoint $x$.
The Yang-Baxter operator constructed from the quandle $Q=x^G$
then leads to the colouring number $\TF{Q} = \F{x}{G}\cdot|Q|$.

We will now move from colouring numbers to colouring polynomials:
consider the quandle extension $\Lambda\curvearrowright\tilde{Q}\to Q$
as defined in \textsection\ref{sub:QuandleCoverings},
where the quandle $Q=Q(G,x)$ is covered by $\tilde{Q}=\tilde{Q}(G,x)$,
and the deck transformation group is $\Lambda = C(x) \cap G'$.
As before, we linearly extend the quandle structure of $\tilde{Q}$ 
to a Yang-Baxter operator $c_{\tilde{Q}}$, and denote the associated 
linear braid group representation by $\rho_{\tilde{Q}}$.
We will, however, not take the total trace 
as before, but rather use the partial trace 
$\tr'\colon \End_\K(\K\tilde{Q}^{\tensor n}) \to \End_\K(\K\tilde{Q})$,
contracting the tensor factors $2,\dots,n$.

\begin{theorem} \label{thm:YBlong}
  Let $(G,x)$ be a finite group such that 
  the conjugacy class $Q=x^G$ generates $G$.
  Let $\tilde{Q} = \tilde{Q}(G,x)$ be the covering quandle
  and let $\rho_{\tilde{Q}}$ be the associated braid group representation.
  Suppose that the knot $K$ is represented by a braid $\beta$.
  Then the partial trace $\tr'( \rho_{\tilde{Q}}(\beta) )\colon \K\tilde{Q} \to \K\tilde{Q}$
  is given by multiplication with $\P{x}{G}(K)$. % $\P{x}{G}([\beta])$.
\end{theorem}

Note that the free left action of $\Lambda$ on $\tilde{Q}$ 
turns $\K\tilde{Q}$ into a free left module over $\K\Lambda$.
In particular, multiplication by $\P{x}{G}(K)$ is a $\K$-linear endomorphism.
If $\K$ is of characteristic $0$, then 
the endomorphism $\tr'(\rho_{\tilde{Q}}(\beta))$ 
uniquely determines $\P{x}{G}(K)$.
% In this sense every knot colouring polynomial 
% can be presented as a Yang-Baxter invariant.

\begin{proof}
  We use the obvious bases $\tilde{Q}$ for $\K\tilde{Q}$ 
  and $\tilde{Q}^n$ for $\K\tilde{Q}^{\tensor{n}}$.
  Each endomorphism $f\colon \K\tilde{Q}^{\tensor n} \to \K\tilde{Q}^{\tensor n}$
  is then represented by a matrix $M^{p_1 p_2\dots p_n}_{q_1 q_2 \dots q_n}$,
  indexed by elements $p_i$ and $q_j$ in the basis $\tilde{Q}$.
  The partial trace $\tr'(f)\colon \K\tilde{Q} \to \K\tilde{Q}$ 
  is given by the matrix $T^{p_1}_{q_1} = \sum M^{p_1 p_2\dots p_n}_{q_1 p_2 \dots p_n}$,
  where the sum is taken over all repeated indices $p_2,\dots,p_n$.
  % By construction, $c_{\tilde{Q}}$ is represented by a permutation matrix,
  % and so is $\rho_{\tilde{Q}}(\beta)$ for every braid $\beta$.
  
  % We first show that only colourings contribute to the trace.
  By construction, each elementary braid $\sigma_i$ 
  acts as a permutation on the basis $\tilde{Q}^n$,
  thus each braid $\beta\in\Br{n}$ is represented 
  by a permutation matrix with respect to this basis.
  We interpret this action as colouring 
  the braid $\beta$ with elements of $\tilde{Q}$: 
  we colour the right ends of the braid 
  with $v = p_1\tensor\cdots\tensor p_n$.
  Moving from right to left, at each crossing the new arc 
  is coloured according to the Wirtinger rule
  as depicted in Figure \ref{fig:ColouredCrossing}.
  We thus arrive at the left ends of the braid being coloured 
  with $\rho(\beta) v = q_1\tensor\cdots\tensor q_n$.
  We conclude that colourings of the braid $\beta$ that satisfy
  the trace conditions $p_2 = q_2, \dots, p_n = q_n$ 
  are in natural bijection with colourings 
  of the corresponding long knot $K$.
  
  We now turn to the remaining indices $p_1$ and $q_1$.
  Let us first consider the special case $p_1 = (x,1)$ 
  and $q_1=(y,\lambda)$. From the preceding argument we see that 
  $T^{p_1}_{q_1}$ % = \sum M^{p_1 p_2\dots p_n}_{q_1 p_2 \dots p_n}$
  equals the number of $\tilde{Q}$-colourings of the long knot $K$ 
  that start with $(x,1)$ and end with $(y,\lambda)$.
  According to Lemma \ref{lem:LiftingLemma},
  such colourings exist only for $y=x$ and $\lambda\in\Lambda$,
  hence we have $q_1 = \lambda \cdot p_1$.
  We conclude that $T^{p_1}_{q_1}$ equals the number 
  of representations $(\pi_{K},m_K,l_K) \to (G,x,\lambda)$.
  In total we get $\tr'(\rho(\beta))\,(p_1) = \P{x}{G}(K) \cdot p_1$.
  
  The preceding construction is equivariant under the right-action
  of the group $G'$ on the covering quandle $\tilde{Q}$.
  According to Lemma \ref{lem:CoveringQuandle} this action 
  is transitive:  for every $p \in \tilde{Q}$ there exists 
  $g \in G'$ and $p = p_1^g$, so we  conclude that 
  $\tr'(\rho(\beta))\,(p) = \P{x}{G}(K) \cdot p$.
  This means that the endomorphism 
  $\tr'(\rho(\beta))\colon \K\tilde{Q} \to \K\tilde{Q}$
  is given by multiplication with $\P{x}{G}(K)$. 
\end{proof}

\begin{remark}
  The partial trace $\tr' \colon \End_\K(\K\tilde{Q}^{\tensor n}) 
  \to \End_\K(\K\tilde{Q})$ corresponds to closing the strands $2,\dots,n$
  of the braid $\beta$, but leaving the first strand open:
  the object thus represented is a long knot.
  The natural setting for such constructions is the category of 
  tangles and its linear representations \cite{Kassel:1995}.
  The previous theorem then says that the long knot $K$ 
  is represented by the endomorphism $\K\tilde{Q} \to \K\tilde{Q}$
  that is given by multiplication with $\P{x}{G}(K)$.

  If we used the complete trace $\tr \colon \End_\K(\K\tilde{Q}^{\tensor n})
  \to \K$ instead, then we would obtain a different invariant 
  $\TF{\tilde{Q}} = \tr\compose\rho_{\tilde{Q}}$.
  By the preceding arguments, $\TF{\tilde{Q}}(K)$ equals $|\tilde{Q}|$ 
  times the number of representations $(\pi_{K},m_K,l_K) \to (G,x,1)$,
  which corresponds to the coefficient of the unit element 
  in the colouring polynomial $\P{x}{G}(K)$.
\end{remark}

\subsection{Colouring polynomials of closed knots} \label{sub:YBclosed}

We will now show how the colouring polynomial $\P{x}{G}$ of closed knots
can be obtained as the trace of a suitable Yang-Baxter representation.
To this end we will modify the construction of the preceding paragraph
in order to replace the partial trace $\tr'$ by the complete trace $\tr$.

We proceed as follows: the quandle $Q=x^G$ admits 
an extension $\Lambda \curvearrowright \tilde{Q} \to Q$ 
as defined in \textsection\ref{sub:QuandleCoverings}.
The quandle structure of $\tilde{Q}$ linearly extends 
to a Yang-Baxter operator $c_{\tilde{Q}}$ on $\K\tilde{Q}$.
The free $\Lambda$-action on $\tilde{Q}$ turns $\K\tilde{Q}$ 
into a free module over $\A = \K\Lambda$. 
If $\Lambda$ is abelian, we can pass to an $\A$-linear operator 
\[
\tilde{c}_Q\colon \K\tilde{Q}\tensor[\A]\K\tilde{Q} \to \K\tilde{Q}\tensor[\A]\K\tilde{Q}
\quad\text{with}\quad \tilde{a} \tensor \tilde{b} 
\mapsto \tilde{b} \tensor ( \tilde{a}\ast\tilde{b} )
\quad\text{for all}\quad \tilde{a}, \tilde{b} \in \tilde{Q}.
\]
The difference between $c_{\tilde{Q}}$ and $\tilde{c}_Q$ 
is that the tensor product is now taken over $\A$, 
which means that everything is bilinear 
with respect to multiplication by $\lambda\in\Lambda$.
In the following theorem and its proof 
all tensor products are to be taken over the ring $\A$,
but for notational simplicity we will write 
$\tensor$ for $\tensor[\A]$.

\begin{theorem} \label{thm:YBclosed}
  If $(G,x)$ is a colouring group such that $\Lambda = C(x)\cap G'$ is abelian,
  then the colouring polynomial $\P{x}{G}\colon \knots\to\Z\Lambda$ is a Yang-Baxter invariant.
  More precisely, the preceding construction yields a Yang-Baxter-Markov operator 
  $(\tilde{c}_Q,\id)$ over the ring $\A = \K\Lambda$, and the associated 
  knot invariant satisfies $\tilde{F}_Q = \varphi\P{x}{G}\cdot|Q|$
  where $\varphi \colon \Z\Lambda \to \K\Lambda$ is the natural 
  ring homomorphism defined by $\varphi(\lambda)=\lambda$
  for all $\lambda \in \Lambda$.
\end{theorem}

If $\K$ is of characteristic $0$, then $\tilde{F}_Q$ is 
equivalent to the knot colouring polynomial $\P{x}{G}$.
If $\K$ is of finite characteristic, then we may lose 
some information and $\tilde{F}_Q$ is usually weaker than $\P{x}{G}$.
In the worst case $|Q|$ vanishes in $\K$ and $\tilde{F}_Q$ becomes trivial.

\begin{proof}
  It is a routine calculation to prove that 
  $\tilde{c}_Q$ is a Yang-Baxter operator over $\A$:
  as before, axiom (Q2) implies that $\tilde{c}_Q$ is an automorphism,
  while axiom (Q3) ensures that $\tilde{c}_Q$ satisfies the Yang-Baxter equation.
  Axiom (Q1) implies the trace condition $\tr_2(\tilde{c}_Q^{\pm1})=\id$,
  hence $(\tilde{c}_Q,\id)$ is a Yang-Baxter-Markov operator. % over $\A$.
  We thus obtain a linear braid group representation 
  $\tilde\rho^n_Q\colon \Br{n} \to \Aut_\A(\K\tilde{Q}^{\tensor n})$,
  whose character $\smash{\tilde{F}_Q} = \tr\compose\tilde{\rho}_Q$
  is Markov invariant and induces a link invariant $\tilde{F}_Q\colon \links\to\A$.
  Restricted to knots we claim that $\tilde{F}_Q = \P{x}{G}\cdot|Q|$.
  The proof of the theorem parallels the proof of Theorem \ref{thm:YBlong},
  but requires some extra care. %It proceeds in four steps.

  % \textit{Representing $\tilde{c}_Q$ by a twisted permutation matrix.}
  To represent $\tilde{c}_Q$ by a matrix, we have to choose a basis
  of $\K\tilde{Q}$ over $\A$. Let $s\colon Q\to\tilde{Q}$ be a section 
  to the central extension $\Lambda\curvearrowright\tilde{Q}\to Q$.
  Then $B=s(Q)$ is a basis of $\K\tilde{Q}$ as an $\A$-module.
  For the basepoint $x$ we can assume $s(x) = (x,1)$, 
  but otherwise there are no canonical choices.
  In general, $s$ will not (and cannot) be a homomorphism of quandles,
  but we have $s(a) \ast s(b) = \lambda(a,b) \cdot s(a\ast b)$
  with a certain $2$-cocycle $\lambda\colon Q\times Q\to \Lambda$ 
  that measures the deviation of $s$ from being a homomorphism.
  % The map $\lambda$ turns out to be a $2$-cocycle in the sense 
  % of quandle cohomology, as discussed in \textsection\ref{sub:QuandleCoverings}.
  Just as $c_Q$ is represented by a permutation matrix, we see that
  $\tilde{c}_Q$ is represented by the same matrix except that the $1$'s 
  are replaced with the elements $\lambda(a,b) \in \Lambda$.
  This is usually called a \emph{monomial matrix} 
  or \emph{generalized permutation matrix}.
  
  % \textit{Identifying eigenvectors and colourings.}
  Since $\K\tilde{Q}$ is a free $\A$-module with finite basis $B = s(Q)$,
  the tensor product $\K\tilde{Q}^{\tensor n}$ is also free and has finite basis $B^n$.
  The trace $\tr\compose\tilde\rho(\beta)$ is calculated 
  as the sum $\sum_{v\in B^n} \langle \tilde\rho(\beta)v | v \rangle$.
  Note that $\tilde\rho(\beta)$ is again a monomial matrix in the sense 
  that each row and each column has exactly one non-zero entry.
  Hence a vector $v\in B^n$ contributes to the trace sum if and only if 
  $\tilde\rho(\beta) v = \lambda(v) v$ with some $\lambda(v)\in\Lambda$.
  It remains to characterize eigenvectors and identify their eigenvalues.
  % The following two steps will characterize eigenvectors
  % and identify their eigenvalues.
  
  Given a braid $\beta\in\Br{n}$ we can interpret the action 
  of $\tilde\rho(\beta)$ as colouring the braid $\beta$: 
  we colour the right ends of the braid with a basis vector $v\in B^n$,
  \[
  v = (a_1,g_1) \tensor (a_2,g_2) \tensor \dots \tensor (a_n,g_n).
  \]
  Moving from right to left, at each crossing 
  the new arc is coloured according to the Wirtinger rule 
  as depicted in Figure \ref{fig:ColouredCrossing}. %  \ref{fig:MeriLongCol}.
  We thus arrive at the left ends of the braid, being coloured with
  \[
  \tilde\rho(\beta) v  = (b_1,h_1) \tensor (b_2,h_2) \tensor \cdots \tensor (b_n,h_n).
  \]
  Since the tensor product is defined over $\A$, 
  we have $\tilde\rho(\beta) v = \lambda(v) v$ 
  if and only if $a_1=b_1, a_2=b_2, \dots, a_n=b_n$.
  Hence each eigenvector $v\in B^n$ naturally corresponds
  to a $Q$-colouring of the closed braid $K=[\beta]$.
  
  % \textit{Identifying eigenvalues and longitudes.}
  In order to identify the eigenvalue $\lambda(v)$,
  we will further assume that $(a_1,g_1) = (x,1)$,
  where $x$ is the basepoint of $G$.
  Such an eigenvector will be called \emph{normalized}.
  % The general case where $(a_1,g_1)$ is arbitrary,
  % can be normalized via conjugation by $g_1^{-1}$.
  % Since everything is equivariant under this action,
  % we obtain the same eigenvalue as in the normalized case.
  Using the tensor product-structure over $\A = \K\Lambda$, we obtain
  \[
  \tilde\rho(\beta) v  
  = (x,\lambda) \tensor (a_2,g_2) \tensor \cdots \tensor (a_n,g_n)
  = \lambda(v) v 
  \]
  as in the proof of Theorem \ref{thm:YBlong}.
  We conclude that each normalized eigenvector $v \in B^n$ 
  with $\tilde\rho(\beta) v = \lambda(v) v$ corresponds to 
  a $\tilde{Q}$-colouring of the long knot, 
  where the first arc is coloured by $(x,1)$
  and the last arc is coloured by $(x,\lambda)$.
  This means that the eigenvalue $\lambda(v)$ 
  is the associated colouring longitude.
  
  % \textit{Calculating the trace via homogeneity.}
  We finally show that $\tilde{F}_Q = \P{x}{G}\cdot|Q|$ by calculating 
  the trace $\sum_{v\in B^n} \langle \tilde\rho(\beta)v | v \rangle$.
  Normalized eigenvectors $v\in \{(x,1)\} \times B^{n-1}$ with 
  $\tilde\rho(\beta) v = \lambda(v) v$ correspond to colourings 
  $\rho \colon (\pi_{K},m_k) \to (G,x)$ with $\rho(l_K) = \lambda(v)$.
  Summing over these vectors only, we thus obtain the 
  colouring polynomial $\P{x}{G}(K)$.
  To calculate the total sum we use again the fact that 
  the right-action of $G'$ on $\tilde{Q}$ is transitive.
  Hence for every $q\in Q$ there exists $g\in G'$ such that $s(q)^g = (x,1)$.
  The action of $g$ induces a bijection between the set of basis vectors 
  $\{s(q)\} \times B^{n-1}$ and $\{(x,1)\} \times B^{n-1}$. 
  Since the preceding trace calculation is $G'$-invariant, each vector
  $v \in \{s(q)\} \times B^{n-1}$ contributes $\P{x}{G}(K)$ to the trace.
  In total we obtain $\tilde{F}_Q = \P{x}{G}\cdot|Q|$, as claimed.
\end{proof}

\subsection{Concluding remarks} \label{sub:Conclusion}

It follows from our construction that $\tilde{c}_Q$ 
is a deformation of the Yang-Baxter operator $c_Q$.
More precisely we have $\tilde{c}_Q(a\tensor b) = \lambda(a,b) \cdot c_Q(a,b)$
for all $a,b\in Q$ with a suitable map $\lambda \colon Q \times Q \to \Lambda$.
Our construction via quandle coverings and central extensions 
provides a geometric interpretation in terms of meridian-longitude information.
This interpretation carries through all steps of our construction,
which finally allows us to interpret the resulting Yang-Baxter 
invariant as a colouring polynomial.

Conversely, it is natural to consider the ansatz 
$\tilde{c}_Q(a\tensor b) = \lambda(a,b) \cdot c_Q(a,b)$ and to ask
which $\lambda$ turn $\tilde{c}_Q$ into a Yang-Baxter operator.
This idea can, though in a restricted form, already 
be found in \cite[Thm.\,4.2.6]{FreydYetter:1989}.
A direct calculation shows that $\tilde{c}_Q$ is 
a Yang-Baxter operator if and only if $\lambda$ 
is a $2$-cocycle in the sense of quandle cohomology.
Moreover, two such deformations will be equivalent
if the cocycles differ by a coboundary. 
% Thus such deformations of $c_Q$ are 
% classified by quandle cohomology $H^2(Q)$.
This observation has been worked out by M.\,Gra\~na \cite{Grana:2002}, 
who independently proved that quandle $2$-cocycle 
state-sum invariants are Yang-Baxter invariants.
% The general theory of Yang-Baxter deformations of $c_Q$ 
% over the power series ring $\K\mathopen{[\![} h \mathclose{]\!]}$ 
% has been developed in \cite{Eisermann:2005}.

%%%%%%%%%%%%%%%%%%%%%%%%%%%%%%%%%%%%%%%%%%%%%%%%%%%%%%%%%%%%%%%%%%%%%%%%%%%%%

\bibliographystyle{amsplain}
\bibliography{colopoly}

%%%%%%%%%%%%%%%%%%%%%%%%%%%%%%%%%%%%%%%%%%%%%%%%%%%%%%%%%%%%%%%%%%%%%%%%%%%%%
\end{document}